\documentclass[11pt]{article}
\usepackage{amsfonts,latexsym,amsmath,amscd,geometry}
\usepackage[utf8]{inputenc}
\usepackage{authblk} 
\usepackage{hyperref} 
\usepackage{amssymb}
\usepackage{latexsym}
\usepackage{xcolor}
\numberwithin{equation}{section}

\newcommand \nc{\newcommand}
\newtheorem{theorem}{Theorem}[section]
\newtheorem{lemma}[theorem]{Lemma}
\newtheorem{proposition}[theorem]{Proposition}

\newtheorem{definition}[theorem]{Definition}

\nc{\ba}{\begin{array}}\nc{\ea}{\end{array}}
\nc{\be}{\begin{eqnarray}}\nc{\ee}{\end{eqnarray}}
\nc{\beq}{\begin{equation}}\nc{\eeq}{\end{equation}}
\nc{\bex}{\begin{eqnarray*}}\nc{\eex}{\end{eqnarray*}}
\nc{\btm}{\begin{theorem}} \nc{\etm}{\end{theorem}}
\nc{\blm}{\begin{lemma}} \nc{\elm}{\end{lemma}}
\nc{\R}{\mathbb{R}}  \nc{\ld}{\lambda}
\nc{\va}{\varphi}
\nc{\ve}{\varepsilon}

\def\pf{\noindent{\bf Proof.\quad}}

\newcommand{\brac}[1]{\left (#1 \right )}
\newcommand \qed {\hfill $\Box$}
\newcommand{\lap}{\Delta }

\newcommand{\ageq}{\succsim}

\renewcommand{\i}{{\bf i}}
\newcommand{\col}[1]{{\color{red} #1}}
\newcommand{\abs}[1]{\left |#1 \right |}

\begin{document}

\title{Partial regularity of heat flow of fractional harmonic maps to spheres}
\author[1]{Hengrong Du \thanks{\href{mailto:hdu@fisk.edu}{hdu@fisk.edu}}}
\author[2]{Tao Huang \thanks{\href{mailto:taohuang@wayne.edu}{taohuang@wayne.edu}}}
\author[3]{Armin Schikorra \thanks{\href{mailto:armin@pitt.edu}{armin@pitt.edu}}}
\author[4]{Changyou Wang \thanks{\href{mailto:wang2482@purdue.edu}{wang2482@purdue.edu}}}

\affil[1]{Department of Mathematics and Computer Science,
Fisk University, Nashville, TN 37208}
\affil[2]{Department of Mathematics, Wayne State University, Detroit, Michigan 48202}
\affil[3]{Department of Mathematics, University of Pittsburgh, Pittsburgh, PA 15260}
\affil[4]{Department of Mathematics, Purdue University, West Lafayette, IN 47907}
\date{}

\maketitle

\begin{abstract} In this paper, we introduce the notation of stationary heat flow of $s$-harmonic maps from $\R^m$ into the unit sphere 
$\mathbb{S}^{\ell-1}\hookrightarrow \mathbb{R}^\ell$, for $0<s<1$,
framed through the parabolic version of Caffarelli-Silvestre's extension. We prove that any stationary heat flow of $s$-harmonic maps
is smooth away from a closed set of  vanishing $(m+2-2s)$-dimensional parabolic Hausdorff measure. 
This extends to the fractional settings the partial regularity
theorem on the heat flow of harmonic maps into spheres by Chen-Li-Lin \cite{CLL} and Feldman \cite{Feldman}.
\end{abstract}

\tableofcontents

\section{Introduction}

Geometric evolution equations play a central role in geometric analysis and the calculus of variations. Among them, the harmonic map heat flow has served as a fundamental paradigm for understanding the interplay between geometry and regularity theory. Given 
two complete smooth Riemannian manifolds $(M,g)$ and $(N,h)\hookrightarrow\R^{\ell}$, the heat flow of harmonic map $u:M\times (0,T)\to N$, $T>0$,  evolves according to
\begin{equation}\label{hf1}
\partial_tu-\Delta u\perp T_uN, 
\end{equation}
and can be interpreted as the (negative) gradient flow of the Dirichlet energy functional
\[
E(u)=\frac12\int_M |\nabla u|^2\,dv_g, \ u: M\to N\hookrightarrow\R^{\ell}.
\]
In their landmark  paper, Eells and Sampson \cite{ES} established the existence of a unique global smooth solution of \eqref{hf1} whenever the target manifold
$(N, h)$ has non-positive curvature. Subsequently, pioneering advances by Struwe \cite{Struwe}, Chen-Struwe \cite{Chen-Struwe}, and Chen-Lin \cite{ChenLin} established the global existence of partially smooth harmonic heat flows
into any compact target manifold $N$ via the development of monotonicity inequalities, local energy inequalities, and $\varepsilon$-priori regularity estimates 
for  Ginzburg-Landau type approximations. Their techniques--most notably the energy monotonicity inequality-laid the foundation for the
partial regularity theories for stationary 
weak solutions of the heat flow of harmonic maps, developed independently by Chen-Li-Lin \cite{CLL} and Feldman \cite{Feldman} for sphere targets, and later
extended by Moser \cite{Moser} for arbitrary compact target manifolds. 

In recent years, nonlocal geometric variational problems driven by fractional operators have attracted considerable interest, including fractional harmonic maps, nonlocal minimal surfaces. These developments naturally motivate the study of geometric evolution equations governed by fractional operators. For $0<s<1$, a natural formulation for the fractional harmonic map flow is
\begin{equation}\label{nonlocalhf1}
\partial_tu+(-\Delta)^su\perp T_uN,
\end{equation}
which corresponds to the negative gradient flow of the fractional Dirichlet energy. While this formulation is natural from both physical and geometric prospectives,
there seems lack of analytic tools to establish its existence and regularity, see Schikorra-Sire-Wang \cite{SSW2017}, Wettstein \cite{Wettstein2022, Wettstein2023},  Sire-Wei-Zheng \cite{SWZ2021},
and Struwe \cite{Struwe2024}  for some related results in this direction.

In this paper, we will investigate an alternative, structurally advantageous formulation of heat flow of fractional harmonic maps, given by
\begin{equation}\label{nonlocalhf2}
(\partial_t-\Delta)^su\perp T_uN
\end{equation}
where $(\partial_t-\Delta)^s$ denotes the fractional power of the parabolic heat operator. 
Both flows \eqref{nonlocalhf1} and \eqref{nonlocalhf2} serve as effective tools for identifying steady states, as both converge to fractional harmonic maps as $t \to \infty$. However, formulation \eqref{nonlocalhf2} offers a distinct analytical advantage: via a parabolic extension, it transforms a purely nonlocal problem on $\mathbb{R}^m \times \mathbb{R}$ into a local, degenerate parabolic equation in one higher spatial dimension, where the degeneracy is dictated by a Muckenhoupt $A_2$-weight and the map is subject to a nonlinear geometric boundary condition. Crucially, the variational structure of the extension framework yields both a Struwe-type energy monotonicity inequality and a local energy estimate. By contrast, the absence of an appropriate monotonicity formula for \eqref{nonlocalhf1} presents a primary obstacle to its  analysis,
especially the issue of partial regularity.

Here, we focus on the fractional harmonic map flow \eqref{nonlocalhf2} targeting the sphere $N = \mathbb{S}^{\ell-1}$. We establish a partial regularity theory for stationary weak solutions by developing analytical tools tailored to the fractional parabolic regime. Our approach combines the semigroup representation of the fractional parabolic operator with its Caffarelli--Silvestre type extension, exploiting both the intrinsic nonlocal structure and the degenerate local PDE in the upper half-space. Under a natural stationarity condition, we establish that both Struwe-type monotonicity formula and local energy estimate hold for the extended flows. We expect that this paradigm—translating boundary non-linearities into extended degenerate parabolic systems—will extend to flows to general target manifolds $N$.

To make this precise, we recall the definition of the fractional heat operator $(\partial_t - \Delta)^s u$ on $\mathbb{R}^m \times \mathbb{R}$, due to Nyström--Sande \cite{NS} and Stinga--Torrea \cite{ST}.
For $0<s<1$ and $u\in \mathcal{D}(H^s)$, we define the Poisson representation formula 
\begin{equation}\label{poisson}
\left(\partial_t-\Delta\right)^s u(x,t)=\iint_{\R^m\times\R_+} (u(x,t)-u(x-z,t-\tau))K_s(z,\tau)\,dyd\tau,
\end{equation}
where the kernel $K_s$ is given by
$$
K_s(z,\tau)=\begin{cases} \frac{1}{(4\pi)^{\frac{m}2}|\Gamma(-s)|} \frac{e^{-\frac{|z|^2}{4\tau}}}{\tau^{\frac{m}2+1+s}}, & \ \forall z\in\R^{m}, \tau>0,\\
0 & \ \ \forall z\in\R^m, \tau\le 0.
\end{cases}
$$ 
and $\Gamma$ denotes the Gamma function. The domain $\mathcal{D}(H^s)$ is defined is
\begin{equation} \label{domain-hs}
\mathcal{D}(H^s)=\Big\{v\in \mathcal{S}'(\R^{m+1}):\  \widehat{v}\in L^1_{\rm{loc}}(\R^{m+1}),  \ (i\tau+|\xi|^2)^{\frac{s}2} \widehat{v}(\xi,\tau)\in L^2(\R^{m+1})\Big\},
\end{equation}
where $\mathcal{S}'(\R^{m+1})$ denotes the space of tempered distributions and $\widehat{v}$ is the Fourier transform of $v$. The domain
$\mathcal{D}(H^s)$ coincides with the anisotropic (homogeneous) Sobolev space 
$$\dot{H}^{s,s/2}_{x,t}(\R^{m+1})
=\Big\{v\in L^1_{\rm{loc}}(\R^{m+1}): \ \big\|v\big\|^2_{\dot{H}^{s,s/2}_{x,t}}\equiv\iint_{\R^{m+1}} (|\tau|^{s}+|\xi|^{2s})|\widehat{v}|^2(\xi,\tau)<+\infty\Big\}.
$$
We define the  the corresponding space for maps into the sphere $\mathbb{S}^{\ell-1}$  by
\[
\dot{H}^{s,s/2}_{x,t}\big(\R^{m+1},\mathbb{S}^{\ell-1}\big)
=\Big\{u=(u^1,\cdots, u^{\ell}):\R^{m+1}\to \mathbb{S}^{\ell-1}\ \big| \ u^i\in \dot{H}^{s,s/2}_{x,t}(\R^{m+1}), \ i=1,\cdots, \ell\Big\}.
\]

We now give the notation of weak solutions to the heat flow of $s$-harmonic maps from $M=\R^m$ into the sphere $\mathbb{S}^{\ell-1}\hookrightarrow\R^{\ell}$.
 
\begin{definition} For $0<s<1$ and an initial data $u_0\in \dot{H}^s(\R^m, \mathbb S^{\ell-1})$,  $u\in \dot{H}^{s,s/2}_{x,t}\big(\R^{m}\times\R,\mathbb{S}^{\ell-1}\big)$ is a {\em weak solution} to  heat flow of $s$-harmonic maps with the initial data $u_0$, 
if it satisfies, in the sense of distributions,
\begin{equation}
\label{eq:fhf1}
\begin{cases}
\displaystyle\left(\partial_t -\Delta\right)^s u =\lambda(x,t) u &\,\,\,\,\hbox{ in }\  \  \R^m \times (0,\infty),\\
u(x,t) = u_0(x)&\,\,\,\,\hbox{ in }\  \  \R^m\times (-\infty,0],
\end{cases}
\end{equation}
where $\lambda$ is the Lagrangian multiplier given by
\[
\lambda(x,t)=\frac12\iint_{\R^m\times\R_+} |u(x,t)-u(x-z,t-\tau)|^2K_s(z,\tau)\,dyd\tau.
\]
Because $\lambda(x,t) u\perp T_u\mathbb{S}^{\ell-1}$,  the first equation of \eqref{eq:fhf1} can be equivalently expressed in the geometric form
\begin{equation}\label{fhf10}
\displaystyle\left(\partial_t -\Delta\right)^s u\perp T_u \mathbb S^{\ell-1}.
\end{equation}
\end{definition}

Following the parabolic Caffarelli--Silvestre extension framework \cite{CaSil, NS, ST},  we reformulate \eqref{eq:fhf1}  into a degenerate parabolic system
in $\R^{m+1}_+\times \R$ with weight $y^a$, $a=1-2s\in (-1,1)$, subject to a nonlinear geometric boundary condition on
$\R^m=\partial \R^{m+1}_+$. Specifically, for $(X,t) = (x,y,t) \in \mathbb{R}^m \times (0,\infty) \times \mathbb{R}$, consider the extension parabolic equation:
\begin{equation}\label{dpe1}
\begin{cases}
y^{a}\partial_t U={\rm{div}}_X(y^aDU) & \ {\rm{in}}\ \ \  \ \R^{m+1}_+\times\R,\\
U(x,0,t)=u(x,t) & \ {\rm{in}}\  \ \ \ \R^m\times \R.
\end{cases}
\end{equation}
For any $u\in \mathcal{D}(H^s)$, there exists a unique smooth solution $U$ of \eqref{dpe1} satisfying the trace condition
\begin{equation}\label{s-po}
-\frac{1}{c_s} y^a \frac{\partial U}{\partial y}\big|_{y=0}=(\partial_t-\Delta)^s u,\ \ {\rm{with}}\ \ c_s=\frac{\Gamma(1-s)}{2^{2s-1}\Gamma(s)}.
\end{equation}
Throughout this paper, we will denote 
\[y^a \frac{\partial U}{\partial y}\Big|_{y=0}=\lim_{y\to 0^+}y^a \frac{\partial U}{\partial y}, \]
where the limit  is understood in the sense of distributions on $\R^m\times \R$. Consequently, the original flow \eqref{eq:fhf1} translates to the local extended system:
\begin{equation}
\label{eq:fhf2}
\begin{cases}
\displaystyle y^a\partial_t U(X,t) = {\rm{div}}_X(y^a DU)(X,t),\,\,\,
&\hbox{ in }\ \ \ \mathbb{R}_+^{m+1}\times \R,\\
U(x,0,t) = u_0(x), \,\,\,\,\,&\hbox{ in }\ \ \ \R^m \times (-\infty,0],\\
\displaystyle y^a\frac{\partial U(X,t)}{\partial y}\Big|_{y=0}
\perp T_u \mathbb S^{\ell-1}, \,\,\,\,\,\,&\hbox{ on }\ \ \ \mathbb{R}^m\times (0,+\infty).
\end{cases}
\end{equation}

By the extension theorem \cite{CaSil} there exists a harmonic extension $\widehat{U}_0\in \dot{H}^1((\R^{m+1}_+, y^a\,dX),\R^{\ell})$ solving
\[
 {\rm{div}}_X(y^a D\widehat{U}_0)=0 \ \ 
\hbox{ in }\ \ \mathbb{R}_+^{m+1}; \  \ \widehat{U}_0(x,0)=u_0(x) \  \ {\rm{on}}\ \ \R^m.
\]
Using $\widehat{U}_0$ as stationary initial data, \eqref{eq:fhf2} (and hence \eqref{eq:fhf1})
becomes the initial and boundary value problem as
\begin{equation}
\label{eq:fhf3}
\begin{cases}
\displaystyle y^a\partial_t U(X,t) = {\rm{div}}_X(y^a DU)(X,t),\,\,\,
&(X,t)\in \mathbb{R}_+^{m+1}\times (0,\infty),\\
U(x,y,0) = \widehat{U}_0(x,y), \,\,\,\,\,& (x,y)\in \R^{m+1}_+,\\
|U|(x,0,t)=1;  \  \displaystyle y^a\frac{\partial U(x, y,t)}{\partial y}\Big|_{y=0}
\perp T_{U(x,0,t)} \mathbb S^{\ell-1}, \,\,\,\,\,\,& (x,t)\in \mathbb{R}^{m}\times (0,+\infty).
\end{cases}
\end{equation}

Throughout this paper, $\mu_a=y^a\,dX$ (and  $\mu_{-a}= y^{-a}\,dX)$)
denotes the weighted Lebesgue measure on $\R^{m+1}_+$.
Hyder-Segatti-Sire-Wang \cite{HSSW} established that for any $u_0\in \dot{H}^s(\R^m,\mathbb{S}^{\ell-1})$, \eqref{eq:fhf3} admits a global weak solution
$$
U\in L^\infty\big(\R_+;  H^1((\R^{m+1}_+,\ \mu_a), \R^{\ell})\big)\cap H^1\big(\R_+;  L^2((\R^{m+1}_+,\ \mu_a), \R^{\ell})\big),
$$
Moreover, when $s=1/2$,
the boundary trace $u=U\big|_{\mathbb{R}^m\times (0,+\infty)}$ of \eqref{eq:fhf1} is smooth outside a closed set of locally finite $(m+1)$-dimensional
parabolic Hausdorff measure,  see also Badman-Cozzi \cite{BC} for a related work.
In a forthcoming paper \cite{GW}, we will extend the main theorems of \cite{HSSW} to all $0<s<1$ with $s\not=1/2$.

We would like to remark that when $s=1/2$,  the heat flow of $1/2$-harmonic map \eqref{eq:fhf1} can be interpreted as the free boundary problem of
the heat flow of harmonic map $U$ from $\mathbb{R}^{m+1}_+$ into $\mathbb{R}^{\ell}$ with the free boundary condition $U(\partial\mathbb{R}^{m+1}_+)\subset 
\mathbb{S}^{\ell-1}$ via \eqref{eq:fhf3}. This view of point has played an important role in the study of steady $1/2$-harmonic map problems,
see Millot-Sire \cite{MS2015}, Millot-Pegon \cite{MP2020}, Mazowiecka-Rodiac-Schikorra \cite{MRS2020}, Millot-Pegon-Schikorra \cite{MPS2021} for relevant results. We would also like to mention the important work
by Da Lio-Rivi\`ere \cite{DR2011} on the regularity for $1/2$-harmonic maps  via commutator estimates.

In this paper, we address the partial regularity  of stationary weak solutions of \eqref{eq:fhf3},
which satisfy a domain-variation stationarity condition analogously to Feldman's framework \cite{Feldman} for the heat flow of harmonic maps.

 Given admissible vector fields $\xi(X,t)\in C^\infty_0\big(\overline{\mathbb R^{m+1}_+}\times (0,\infty), \mathbb R^{m+1}\big)$, with $\xi^{m+1}\ge 0$,
and test functions $0\le \theta (X,t)\in C^\infty_0\big(\overline{\mathbb R^{m+1}_+}\times (0,\infty)\big)$,
 the inner variation of $U^\tau$ for $\tau>0$ is defined as
$$
U^\tau(X,t)=U\big(X+\tau\xi(X, t), t+\tau\theta(X,t)\big), \ (X,t)\in\overline{\mathbb R^{m+1}_+}\times [0,\infty).
$$

\begin{definition} \label{stationary00}  A map $U:\overline{\mathbb R^{m+1}_+}\times (0,\infty)\to \mathbb R^{\ell}$,
belonging to $L^\infty\big(\R_+;  \dot{H}^1((\R^{m+1}_+; \mu_a), \R^{\ell})\big)\cap \dot{H}^1\big(\R_+; L^2((\R^{m+1}_+; \mu_a), \R^{\ell})\big)$,
is called a stationary solution of the heat flow of $s$-harmonic maps into $\mathbb{S}^{\ell-1}$, if
$U$ is not only a weak solution of \eqref{eq:fhf3} but also satisfies the following variational inequality:
\begin{equation}\label{stationary000}
2\iint_{\mathbb R^{m+1}_+\times \mathbb R_+} \partial_tU \cdot\frac{d}{d\tau}\big|_{\tau=0} U^\tau \mu_adt
+\frac{d}{d\tau}\big|_{\tau=0}\iint_{\mathbb R^{m+1}_+\times\mathbb R_+} \big|DU^\tau\big|^2(X,t)\,\mu_adt\le 0
\end{equation}
for all admissible $\xi$, with $\xi^{m+1}\ge 0$ and $\xi^{m+1}\big|_{y=0}=0$, and $0\le\theta\in C^\infty_0\big(\overline{\mathbb R^{m+1}_+}\times (0,\infty)\big)$.
\end{definition} 

Stationary solutions satisfy two pivotal structural properties: 
\begin{enumerate}
\item [(i)] Struwe's type monotonicity inequality for $y^a$-weighted energy;
and 
\item [(ii)] Local $y^a$-weighted energy inequality. 
\end{enumerate}
These inequalities serve as the primary source for proving 
partial regularity. Our main theorem is as follows.

\begin{theorem} \label{partial_reg} For $0<s<1$, let $U\in L^\infty\big(\R_+;  \dot{H}^1((\R^{m+1}_+; \mu^a), \R^{\ell})\big)
\cap \dot{H}^1\big(\R_+; L^2((\R^{m+1}_+; \mu_a), \R^{\ell})\big)$
be a stationary solution of the heat flow of harmonic maps into $\mathbb{S}^{\ell-1}$ in the sense of Definition~\ref{stationary00}.. 
Then there exists a closed set $\Sigma\subset \mathbb R^m\times (0,\infty)$, with $\mathcal{H}^{m+2-2s}(\Sigma)=0$, such that 
\[
u=U\big|_{\partial\R^{m+1}_+\times (0,\infty)}\in C^\infty(\mathbb R^m\times (0,\infty)\setminus \Sigma,\mathbb S^{\ell-1}).
\]
\end{theorem}
Here $\mathcal{H}^{m+2-2s}$ denotes $(m+2-2s)$-dimensional Hausdorff measure with respect to the parabolic distance:
\[
d\big((x,t), (y,s)\big)=\max\Big\{|x-y|, \ \sqrt{|t-s|}\Big\},\ \  \forall (x,t), (y,s) \in \R^{m}\times (0,\infty).
\]

Theorem \ref{partial_reg} is deduced via an $\varepsilon_0$-regularity criterion under the smallness of the weighted Morrey norm $\|\partial_t U\|_{M^{2,4},\mu_a dt} + \|DU\|_{M^{2,2},\mu_a dt}$. Beyond monotonicity inequality and 
local energy inequality, a crucial structural ingredient is the hidden algebraic structure of \eqref{eq:fhf3}, which can be rewritten
as (see Lemma \ref{alg0} below)
\begin{align}
y^a\partial_t U^k-{\rm{div}}_X(y^a DU^k)&=\frac12 y^a\partial_t(|U|^2-1)U^k +y^a DU^i\cdot\Omega^{ki}-\frac12 {\rm{div}}_X(y^aD(|U|^2-1)U^k) \nonumber\\
&\qquad\qquad+y^aD(|U|^2-1)\cdot DU^k \  \ \ {\rm{in}}\ \ \ \R^{m+1}_+\times (0,\infty), \label{eq:fhf4}\\
y^a \frac{\partial U^k}{\partial y}\Big|_{y=0}&=\frac12 y^a\frac{\partial}{\partial y}(|U|^2-1) U^k\Big|_{y=0} 
\ \ \ \ \ \ \ \ \ \  {\rm{on}}\ \ \ \partial \R^{m+1}_+\times (0,\infty). \label{bdy_cond}
\end{align}
Here the matrix-valued connection $1$-form $\Omega^{ij}=U^i DU^j-U^j DU^i$ satisfies 
\[
{\rm{div}}_X(y^a\Omega^{ij})=y^a \big(U^i\partial_t U^j-U^j\partial_t U^i\big).
\]
This identity on divergence implies that $y^a DU^i\cdot\Omega^{ki}$ exhibits a div-curl structure so that it is embedding into the (local) Hardy space in $\R^{m+1}_+\times (0,\infty)$. Exploiting Hardy--BMO duality yields an improved energy estimate of $U$ (Theorem \ref{hybrid1}), forming a key step in 
establishing Theorem \ref{partial_reg}.

\subsection*{Organization of the Paper}
The paper is organized as follows. In Section 2, we rewrite system \eqref{eq:fhf3} to uncover the hidden algebraic structure. In Section 3, we establish the improved energy estimates via Hardy--BMO duality. In Section 4, we derive the $y^a$-weighted energy monotonicity and local $y^a$-weighted energy inequalities under the stationarity condition. In Section 5, we execute the blow-up analysis and establish the initial Hölder continuity ($C^{\alpha_s}$). In Section 6, we perform the higher-order bootstrap argument to establish $C^\infty$ regularity off the singular set. Finally, Section 7 (Appendix) collects some key properties of weighted function spaces and integral estimates.

\section{Structure of Nonlinearities of the Heat Flow of $s$-Harmonic Maps}

In this section, we provide the hidden algebraic structure of system \eqref{eq:fhf3} governing the heat flow of $s$-harmonic maps. This structural property will play a critical role in deriving the Hardy--BMO duality estimates and establishing the main partial regularity theorem. Similar algebraic structures were previously identified by Mazowiecka-Rodiac-Schikorra \cite{MRS2020} in the context of free boundary problems for $p$-harmonic maps.

\begin{lemma}\label{lem:halfflow1}
Let $U: \mathbb{R}^{m+1}_+ \times \mathbb{R}_+ \to \mathbb{R}^\ell$ be a weak solution to \eqref{eq:fhf3} with boundary trace 
$u = U\big|_{\partial\mathbb{R}^{m+1}_+ \times \mathbb{R}_+}:  {\partial\mathbb{R}^{m+1}_+ \times \mathbb{R}_+}\to \mathbb{S}^{\ell-1}$. Then for any test function $\Phi \in C^\infty_0(\overline{\mathbb{R}^{m+1}_+} \times \mathbb{R}_+)$ and each $k=1, \cdots, \ell$, we have
\begin{align}\label{eq:halfflow1}
&\iint_{\mathbb{R}^{m+1}_+ \times \mathbb{R}_+} \left( \partial_t U^k \Phi + DU^k \cdot D\Phi \right) \mu_a \mathrm{d}t \nonumber\\
& \qquad\qquad\qquad\qquad\qquad\qquad= \iint_{\mathbb{R}^{m+1}_+ \times \mathbb{R}_+} \left( \partial_t U^i (U^k U^i \Phi) + DU^i \cdot D(U^k U^i \Phi) \right) \mu_a \mathrm{d}t.
\end{align}
\end{lemma}

\pf
Applying the equation $y^a \partial_t U = \mathrm{div}_X (y^a DU)$ in $\mathbb{R}^{m+1}_+ \times \mathbb{R}_+$ and integrating by parts, the difference between the two sides of \eqref{eq:halfflow1} can be written as
\begin{align*}
&\iint_{\mathbb{R}^{m+1}_+ \times \mathbb{R}_+} \left( \partial_t U^k \Phi + DU^k \cdot D\Phi \right) \mu_a \mathrm{d}t - \iint_{\mathbb{R}^{m+1}_+ \times \mathbb{R}_+} \left( \partial_t U^i (U^k U^i \Phi) + DU^i \cdot D(U^k U^i \Phi) \right) \mu_a \mathrm{d}t \\
&\qquad = \iint_{\mathbb{R}^{m+1}_+ \times \mathbb{R}_+} y^a \partial_t U^i \left( \delta_{ki} - U^k U^i \right) \Phi \, \mathrm{d}X \mathrm{d}t - \iint_{\mathbb{R}^{m+1}_+ \times \mathbb{R}_+} y^a DU^i \cdot D \Big( (\delta_{ki} - U^k U^i) \Phi \Big) \mathrm{d}X \mathrm{d}t \\
&\qquad = -\iint_{\partial \mathbb{R}^{m+1}_+ \times \mathbb{R}_+} \lim_{y \to 0^+} \left( y^a \frac{\partial U^i}{\partial y} \right) \left( \delta_{ki} - U^k U^i \right) \Phi \, \mathrm{d}x \mathrm{d}t.
\end{align*}
Since $|u(x,t)| = |U(x,0,t)| = 1$ on $\partial \mathbb{R}^{m+1}_+ \times \mathbb{R}_+$, contractive sum over $i$ gives
\[
u^i (\delta_{ki} - u^k u^i) = u^k - |u|^2 u^k = 0 \quad \text{on } \partial \mathbb{R}^{m+1}_+ \times \mathbb{R}_+.
\]
Consequently, the vector field $\big( \delta_{k1} - u^k u^1, \, \dots, \, \delta_{k\ell} - u^k u^\ell \big)$ is tangent to the sphere $\mathbb{S}^{\ell-1}$ at $u(x,t)$. Recalling the geometric Neumann boundary condition in \eqref{eq:fhf3}, namely
\[
\lim_{y \to 0^+} y^a \frac{\partial U}{\partial y} \perp T_u \mathbb{S}^{\ell-1},
\]
it follows immediately that the boundary integral vanishes:
\[
\iint_{\partial \mathbb{R}^{m+1}_+ \times \mathbb{R}_+} \lim_{y \to 0^+} \left( y^a \frac{\partial U^i}{\partial y} \right) \left( \delta_{ki} - U^k U^i \right) \Phi \, \mathrm{d}x \mathrm{d}t = 0.
\]
This completes the proof of \eqref{eq:halfflow1}.
\qed

\medskip

Recall the matrix-valued connection $1$-form $(\Omega^{ij})_{1 \le i, j \le \ell}$ on $\overline{\mathbb{R}^{m+1}_+} \times \mathbb{R}_+$ defined in Section 1:
\[
\Omega^{ij} = U^i DU^j - U^j DU^i.
\]
By construction, $(\Omega^{ij})$ is skew-symmetric, i.e., $\Omega^{ij} = -\Omega^{ji}$. Its conservation-law properties are established in the following lemma:

\begin{lemma}\label{alg0}
Let $U: \mathbb{R}^{m+1}_+ \times \mathbb{R}_+ \to \mathbb{R}^\ell$ be a weak solution to \eqref{eq:fhf3}. Then for any $\Phi \in C^\infty_0(\overline{\mathbb{R}^{m+1}_+} \times \mathbb{R}_+)$ and each $1 \le i, j \le \ell$, we have
\begin{equation}\label{alg1}
\iint_{\mathbb{R}^{m+1}_+ \times \mathbb{R}_+} \Omega^{ij} \cdot D\Phi \, \mu_a \mathrm{d}t = -\iint_{\mathbb{R}^{m+1}_+ \times \mathbb{R}_+} \left( U^i \partial_t U^j - U^j \partial_t U^i \right) \Phi \, \mu_a \mathrm{d}t.
\end{equation}
\end{lemma}

\pf
Integrating by parts on the left-hand side of \eqref{alg1} yields
\begin{align*}
\iint_{\mathbb{R}^{m+1}_+ \times \mathbb{R}_+} \Omega^{ij} \cdot D\Phi \, \mu_a \mathrm{d}t 
&= -\iint_{\mathbb{R}^{m+1}_+ \times \mathbb{R}_+} \mathrm{div}_X (y^a \Omega^{ij}) \Phi \, \mathrm{d}X \mathrm{d}t \\
&\quad - \iint_{\partial\mathbb{R}^{m+1}_+ \times \mathbb{R}_+} \lim_{y \to 0^+} y^a \left( U^i \frac{\partial U^j}{\partial y} - U^j \frac{\partial U^i}{\partial y} \right) \Phi \, \mathrm{d}x \mathrm{d}t.
\end{align*}
Using the equation $y^a \partial_t U = \mathrm{div}_X (y^a DU)$ satisfied in the interior $\mathbb{R}^{m+1}_+ \times \mathbb{R}_+$, a direct computation shows that
\[
\mathrm{div}_X (y^a \Omega^{ij}) = U^i \mathrm{div}_X (y^a DU^j) - U^j \mathrm{div}_X (y^a DU^i) = y^a \left( U^i \partial_t U^j - U^j \partial_t U^i \right).
\]
For the boundary integral, since $u = U\big|_{\partial\mathbb{R}^{m+1}_+} \in \mathbb{S}^{\ell-1}$ and $\lim_{y \to 0^+} y^a \frac{\partial U}{\partial y} \perp T_u \mathbb{S}^{\ell-1}$, the normal derivative at the boundary is parallel to $u(x,t)$:
\[
\lim_{y \to 0^+} y^a \frac{\partial U}{\partial y} = \left\langle \lim_{y \to 0^+} y^a \frac{\partial U}{\partial y}, \, u \right\rangle u \quad \text{on } \partial\mathbb{R}^{m+1}_+ \times \mathbb{R}_+.
\]
Therefore, componentwise on the boundary $\partial\mathbb{R}^{m+1}_+ \times \mathbb{R}_+$, we have
\[
u^i \lim_{y \to 0^+} y^a \frac{\partial U^j}{\partial y} - u^j \lim_{y \to 0^+} y^a \frac{\partial U^i}{\partial y} 
= u^i u^j \left\langle \lim_{y \to 0^+} y^a \frac{\partial U}{\partial y}, \, u \right\rangle - u^j u^i \left\langle \lim_{y \to 0^+} y^a \frac{\partial U}{\partial y}, \, u \right\rangle = 0.
\]
Hence the boundary contribution vanishes, completing the proof of \eqref{alg1}.
\qed

\medskip

Combining Lemmas \ref{lem:halfflow1} and \ref{alg0}, we obtain the main algebraic identity for the $s$-harmonic map heat flow.

\begin{lemma}\label{lem:halfflow2}
Let $U: \mathbb{R}^{m+1}_+ \times \mathbb{R}_+ \to \mathbb{R}^\ell$ be a weak solution to \eqref{eq:fhf3}. Then for any $\Phi \in C^\infty_0(\overline{\mathbb{R}^{m+1}_+} \times \mathbb{R}_+)$ and $1 \le k \le \ell$, we have
\begin{align}\label{halfflow2}
\iint_{\mathbb{R}^{m+1}_+ \times \mathbb{R}_+} \left( \partial_t U^k \Phi + DU^k \cdot D\Phi \right) \mu_a \mathrm{d}t
&= \frac{1}{2} \iint_{\mathbb{R}^{m+1}_+ \times \mathbb{R}_+} \partial_t (|U|^2 - 1) U^k \Phi \, \mu_a \mathrm{d}t \nonumber \\
&\quad + \iint_{\mathbb{R}^{m+1}_+ \times \mathbb{R}_+} DU^i \cdot \Omega^{ki} \Phi \, \mu_a \mathrm{d}t \nonumber \\
&\quad + \frac{1}{2} \iint_{\mathbb{R}^{m+1}_+ \times \mathbb{R}_+} D(|U|^2 - 1) U^k \cdot D\Phi \, \mu_a \mathrm{d}t \nonumber \\
&\quad + \iint_{\mathbb{R}^{m+1}_+ \times \mathbb{R}_+} D(|U|^2 - 1) \cdot DU^k \Phi \, \mu_a \mathrm{d}t.
\end{align}
\end{lemma}

\pf
We expand the integrands on the right-hand side of identity \eqref{eq:halfflow1} in Lemma \ref{lem:halfflow1}. First, for the temporal derivative term, we write:
\[
\partial_t U^i (U^k U^i \Phi) = \left( U^i \partial_t U^i \right) U^k \Phi = \frac{1}{2} \partial_t (|U|^2 - 1) U^k \Phi.
\]
Next, expanding the spatial gradient term yields:
\begin{align*}
DU^i \cdot D(U^k U^i \Phi) 
&= DU^i \cdot \left( U^i U^k D\Phi + U^i \Phi DU^k + U^k \Phi DU^i \right) \\
&= \left( U^i DU^i \right) \cdot D\Phi \, U^k + \left( U^i DU^i \right) \cdot DU^k \, \Phi + |DU|^2 U^k \Phi.
\end{align*}
Using $U^i DU^i = \frac{1}{2} D(|U|^2 - 1)$, the first two terms become $\frac{1}{2} D(|U|^2 - 1) U^k \cdot D\Phi$ and $\frac{1}{2} D(|U|^2 - 1) \cdot DU^k \Phi$, respectively. 

For the third term, we re-express $|DU|^2 U^k \Phi$ using the connection $1$-form $\Omega^{ki} = U^k DU^i - U^i DU^k$:
\begin{align*}
DU^i \cdot \Omega^{ki} 
&= DU^i \cdot \left( U^k DU^i - U^i DU^k \right) \\
&= |DU|^2 U^k - \left( U^i DU^i \right) \cdot DU^k \\
&= |DU|^2 U^k - \frac{1}{2} D(|U|^2 - 1) \cdot DU^k.
\end{align*}
Rearranging this identity gives
\[
|DU|^2 U^k \Phi = DU^i \cdot \Omega^{ki} \Phi + \frac{1}{2} D(|U|^2 - 1) \cdot DU^k \Phi.
\]
Summing all the expanded pieces, the right-hand side of \eqref{eq:halfflow1} transforms into:
\begin{align*}
&\iint_{\mathbb{R}^{m+1}_+ \times \mathbb{R}_+} \left( \partial_t U^i (U^k U^i \Phi) + DU^i \cdot D(U^k U^i \Phi) \right) \mu_a \mathrm{d}t \\
&\qquad = \frac{1}{2} \iint_{\mathbb{R}^{m+1}_+ \times \mathbb{R}_+} \partial_t (|U|^2 - 1) U^k \Phi \, \mu_a \mathrm{d}t + \iint_{\mathbb{R}^{m+1}_+ \times \mathbb{R}_+} DU^i \cdot \Omega^{ki} \Phi \, \mu_a \mathrm{d}t \\
&\qquad\quad + \frac{1}{2} \iint_{\mathbb{R}^{m+1}_+ \times \mathbb{R}_+} D(|U|^2 - 1) U^k \cdot D\Phi \, \mu_a \mathrm{d}t + \iint_{\mathbb{R}^{m+1}_+ \times \mathbb{R}_+} D(|U|^2 - 1) \cdot DU^k \Phi \, \mu_a \mathrm{d}t.
\end{align*}
Substituting this back into \eqref{eq:halfflow1} establishes \eqref{halfflow2}, completing the proof.
\qed

  \section{Hybrid inequality under small $y^a$-weighted Morrey norm}

In this section, we will derive a hybrid energy inequality for the heat flow of $s$-harmonic map under the smallness condition of $y^a$-weighted 
Morrey norm. First, we recall the definition of $y^a$-weighted Morrey spaces in the parabolic upper half-space.

\begin{definition} Let $O\subset\mathbb R^{m+1}\times (0,\infty)$ be an open set, and
denote $O_+=O\cap \big(\overline{\mathbb R^{m+1}_+}\times (0,\infty)\big)$. For $0<\lambda<m+3+a$, 
let $M^{2,\lambda}(O_+,\ \mu_adt)$ be the space consisting of all measurable functions $f\in L^1_{\rm{loc}}(O_+,\mu_adt)$ such that 
\begin{align}\label{morrey}
\big\|f\big\|_{M^{2,\lambda}(O_+, \  \mu_a dt)}^2=\sup_{(x,t)\in O_+, r>0} \Big\{r^{\lambda-(m+3+a)}\iint_{P_r(x,t)\cap O_+} |f|^2\mu_a dt\Big\}<\infty.
\end{align}
\end{definition}
Throughout this paper,  denote by $P_r(X,t)=B_r(X)\times (t-r^2, t+r^2)$ the parabolic ball centered at $(X,t)$ with radius $r$, 
and $P_r^\pm(X,t)=P_r(X,t)\cap (\R^{m+1}_\pm\times \R)$. The weighted average of $U$ over $P_r^+(z_0)$ with respect to $\mu_adt$
is defined by 
$$(U)^a_{z_0,r,+}=\displaystyle\frac{1}{|{P_{r}^+(z_0)}|_a}\iint_{P_{r}^+(z_0)} U\mu_adt,
\ {\rm{where}}\ |{P_{r}^+(z_0)}|_a=\displaystyle\iint_{P_{r}^+(z_0)}\mu_adt.$$
When $a=0$, $|P_r^+(z_0)|_0$ and $(U)^0_{z_0,r+}$ reduce to the unweighted the Lebesgue measure $|P_r^+(z_0)|$ and the standard
 average $U_{z_0, r,+}$ respectively.

Before proving the main hybrid inequality, we establish a crucial point-wise distance estimate showing that small Morrey bounds force $U$ to stay close to the target manifold $\mathbb{S}^{\ell-1}$.

\begin{lemma} \label{dist-est0} For any $z_0=(x_0, 0, t_0)\in\partial\mathbb{R}^{m+1}_+\times (0,\infty)$ and $0<r_0<\sqrt{t_0}$, we have 
\begin{align}\label{dist-est}
\sup_{z\in P_{\frac{r_0}4}^+(z_0)}{\rm{dist}}\big(U(z),\mathbb{S}^{\ell-1}\big) 
\le C\Big[\big\|DU\big\|_{M^{2,2}(P^+_{r_0}(z_0),\ \mu_adt)}+\big\|\partial_t U\big\|_{M^{2,4}(P^+_{r_0}(z_0),\ \mu^adt)}\Big].
\end{align}
\end{lemma}
\pf Pick any  $z=(x,y,t)\in P_{\frac{r_0}4}^+(z_0)$. Since $U(x,0, t)\in \mathbb S^{\ell-1}$, we may assume $y>0$. Since $U$ solves the heat equation in 
the interior of $\mathbb R^{m+1}_+$
and $P_{\frac{y}2}(z)\Subset P_{y}(z)\subset P^+_{\frac{r_0}2}(z_0)$, it follows from
the interior regularity estimates for parabolic equations that
\begin{align}\label{poin1}
|U(z)-(U)_{z, \frac{y}2}|&\le C\Big(|y|^{2-(m+3+a)}\iint_{P_{|y|}(z)}(|DU|^2+|y|^2||\partial_t U|^2)\mu_adt\Big)^\frac12\nonumber\\
&\le C\Big(\big\|DU\big\|_{M^{2,2}(P^+_{r_0}(z_0),\ \mu_adt)}+\big\|\partial_t U\big\|_{M^{2,4}(P^+_{r_0}(z_0),\ \mu_adt)}\Big).
\end{align}
Set $\bar{z}=(x,0,t)$. By inclusion, $P_{\frac{y}2}(z)\subset P^+_{\frac{3y}2(\bar{z})}\subset P_{r_0}^+(z_0)$.
Applying H\"older's inequality,  and the fact that $y^a$ is a $A_2$-weight, and the $y^a$-weighted Poincar\'e inequality, we obtain 
\begin{align}\label{poin2}
\big|(U)_{z, \frac{y}2}-(U)_{\bar{z}, \frac{3y}2, +}\big|
&\le \frac{3^{m+2}}{|P^+_{\frac{3y}2}(\bar{z})|}\iint_{P^+_{\frac{3y}2}(\bar{z})}|U-(U)_{\bar{z}, \frac{3y}2, +}|\,dXdt\nonumber\\
&\le {C}\Big(\frac{1}{\displaystyle |P^+_{\frac{3y}2}(\bar{z})|_a}\iint_{P^+_{\frac{3y}2}(\bar{z})} |U-(U)^a_{\bar{z}, \frac{3y}2, +}|^2\mu_adt\Big)^\frac12\nonumber\\
&\le C \Big(\big\|DU\big\|_{M^{2,2}(P^+_{r_0}(z_0),\ \mu_adt)}+\big\|\partial_t U\big\|_{M^{2,4}(P^+_{r_0}(z_0),\ \mu_adt)}\Big).
\end{align} 
Since $u(x,t)=U(x,0,t)\in\mathbb{S}^{\ell-1}$, for any boundary point $w \in  \partial \mathbb{R}^{m+1}_+\times \R_+$ we have
\begin{align*}
{\rm{dist}}\Big((U)_{\bar{z}, \frac{3y}2, +}, \mathbb{S}^{\ell-1}\Big)\le \Big|(U)_{\bar{z}, \frac{3y}2, +}-u(w)\Big|
\le Cy^{-(m+3)}\iint_{P^+_{\frac{3y}2}(\bar{z})}\Big|U-u(w)|
\end{align*}
Averaging $w$ over $\widehat{P}_{\frac{3y}2}(\bar z):=D_{\frac{3y}2}(x)\times [t-\frac{9y^2}4, t+\frac{9y^2}4]$, the parabolic ball in $\R^m\times \R$, 
applying the triangle inequality, the weighted Poincar\'e inequality \eqref{weight-poincare}, and the trace inequality \eqref{poincare2}, we find
\begin{align}\label{poin3}
{\rm{dist}}\Big((U)_{\bar{z}, \frac{3y}2, +}, \mathbb{S}^{\ell-1}\Big)
&\le Cy^{-(m+3)}\iint_{P^+_{\frac{3y}2}(\bar{z})}\Big|U-(u)_{\widehat{P}_{\frac{3y}2}(\bar{z})}\Big|\nonumber\\
&+Cy^{-(m+2)}\iint_{\widehat{P}_{\frac{3y}2}(\bar z)}
\Big|u-(u)_{\widehat{P}_{\frac{3y}2}(\bar{z})}\Big|\nonumber\\
&\le C\Big(\big\|DU\big\|_{M^{2,2}(P^+_{r_0}(z_0),\ \mu_adt)}+\big\|\partial_t U\big\|_{M^{2,4}(P^+_{r_0}(z_0),\ \mu_adt)}\Big).
\end{align}
Combining \eqref{poin1},  \eqref{poin2}, and \eqref{poin3} completes the proof of \eqref{dist-est}
\qed

\medskip

Equipped with Lemma \ref{dist-est0}, we establish the key hybrid inequality under a smallness condition.
\begin{theorem} \label{hybrid1}There exist $\varepsilon_0>0$ and $\theta_0>0$ such that for any $z_0=(x_0,t_0)\in\partial\mathbb{R}^{m+1}_+\times (0,\infty)$, if $U\in L^\infty_tH^1_X(P_{r_0}^+(z_0), \mu_adt)$ satisfies
\begin{equation}\label{small1}
r_0^{-(m+1+a)}\iint_{P_{r_0}^+(z_0)}|DU|^2\mu_adt
\le\varepsilon_0^2,
\end{equation}
 then for all $0<r\le \theta_0r_0$, the following hybrid energy inequality 
\begin{align}\label{hybrid2}
&r^{-(m+1+a)}\iint_{P^+_r(z_0)}  |DU|^2\mu_adt\nonumber\\
&\qquad\qquad\qquad\le C\varepsilon_0 r^{-(m+1+a)} \iint_{P^+_{2r}(z_0)} (r^2|\partial_t U|^2+|DU|^2)\mu_adt\nonumber\\
&\qquad\qquad\qquad\ \ +C\varepsilon_0 r^{-(m-1+a)}\sup_{t\in [t_0-4r^2, t_0+4r^2]}\int_{B_{2r}^+(x_0)} |DU|^2\mu_a\nonumber\\ 
&\qquad\qquad\qquad\ \ +C\varepsilon_0^{-1} r^{-(m+3+a)}\iint_{P^+_{2r}(z_0)} \Big(|U-(U)^a_{z_0,2r,+}|^2+ |U-(U)^a_{z_0,2r,+}|^4\Big)\mu_adt\nonumber\\
&\qquad\qquad\qquad\ \ +C\varepsilon_0^{-1}r^{-(m+3+a)} \iint_{P^+_{2r}(z_0)} (|U|^2-1)^2\mu_adt,
\end{align}
holds for some $C=C(m)>0$. 
\end{theorem}  
\pf Without loss of generality, assume $z_0=(0,1)$ and $r_0=\frac12$. For $0<r<\frac18$, let $\eta\in C^\infty_0(P_1(0,1))$ be a  parabolic
cut off function for $P_r(0,1)$ satisfying $0\le \eta\le 1$,
$\eta\equiv 1$ in $P_r(0,1)$, $\eta\equiv 0$ outside $P_{2r}(0,1)$, and $|D^2\eta|+|D\eta|^2+|\partial_t \eta|\le 16r^{-2}$. 

Substituting test function $\Phi=(U^k-(U^k)^a_{(0,1), 2r, +})\eta^2$ into 
\eqref{halfflow2}, we obtain
\begin{align*}
\iint_{P_1^+(0,1)}|DU|^2\eta^2\mu_adt
&=-\iint_{P_1^+(0,1)}\partial_tU^k (U^k-(U^k)^a_{(0,1), 2r, +})\eta^2\mu_adt\\
&\ \ \ -\iint_{P_1^+(0,1)}DU^k (U^k-(U^k)^a_{(0,1), 2r, +})\cdot D\eta^2\mu_adt\\
&\ \ \ +\iint_{P_1^+(0,1)} \frac12 \partial_t(|U|^2-1) U^k(U^k-(U^k)^a_{(0,1), 2r, +})\eta^2\mu_adt\\
&\ \ \ +\iint_{P_1^+(0,1)} \frac12 D(|U|^2-1) U^k \cdot D((U^k-(U^k)^a_{(0,1), 2r, +})\eta^2)\mu_adt\\
&\ \ \ +\iint_{P_1^+(0,1)} D(|U|^2-1)\cdot DU^k(U^k-(U^k)^a_{(0,1), 2r, +})\eta^2\mu_adt\\
&\ \ \ +\iint_{P_1^+(0,1)} DU^i\cdot\Omega^{ki}(U^k-(U^k)^a_{(0,1), 2r, +})\eta^2\mu_adt\\
&=I+II+III+IV+V+VI.
\end{align*} 
Before estimating  $I, II, III, IV, V$, recall Lemma 3.2, \eqref{small1}, Lemma \ref{local_monotonicity},  and Lemma \ref{local_energy_ineq}  imply
that there exists $\theta_0\in (0,\frac14)$ such that
\begin{equation}\label{dist-est1}
\big|1-|U(z)|\big|={\rm{dist}}(U(z), \mathbb{S}^{\ell-1})\le C\varepsilon_0, \ \forall z\in P^+_{2\theta_0}(0,1).
\end{equation} 
In particular, $|U|\le 2$ on $P^+_{2\theta_0}(0,1)$ if $\varepsilon_0>0$ is chosen to be sufficiently small.

For $0<r<\min\{\frac18, \frac{\theta_0}2\}$, by Young's inequality, the terms $I, II, III$ are bounded as follows
\begin{align*}
\big|I\big|\le \varepsilon_0 (2r)^2\iint_{P^+_{2r}(0,1)}|\partial_t U|^2\mu_adt
+C\varepsilon_0^{-1} r^{-2}\iint_{P^+_{2r}(0,1)} \big|U-(U)^a_{(0,1), 2r, +}\big|^2\mu_adt,
\end{align*}
\begin{align*}
\big|II\big|\le \varepsilon_0 \iint_{P^+_{2r}(0,1)} |DU|^2\mu_adt+C\varepsilon_0 r^{-2}\iint_{P^+_{2r}(0,1)} \big|U-(U)^a_{(0,1), 2r, +}\big|^2\mu_adt,
\end{align*}
\begin{align*}
\big|III\big|&\le \frac14\varepsilon_0(2r)^2\|U\|_{L^\infty(P^+_{\frac14}(0,1))}^2 \iint_{P^+_{2r}(0,1)} |\partial_t U|^2\mu_adt
+C\varepsilon_0 r^{-2}\iint_{P^+_{2r}(0,1)} \big|U-(U)^a_{(0,1), 2r, +}\big|^2\mu_adt\\
&\le \varepsilon_0 (2r)^2\iint_{P^+_{2r}(0,1)} |\partial_t U|^2\mu_adt+C\varepsilon_0 r^{-2}\iint_{P^+_{2r}(0,1)} \big|U-(U)^a_{(0,1), 2r, +}\big|^2\mu_adt.
\end{align*}
To estimate $IV$, observe that $(|U|^2-1)\eta^2=0$ on $\partial P^+_1(0,1)$. Integrating by parts yields
\begin{align*}
&IV=\iint_{P_1^+(0,1)} \big[D(|U|^2-1) \langle U, DU\rangle\eta^2+D(|U|^2-1)\langle U, U-(U)^a_{(0,1), 2r, +}\rangle D\eta^2\big]\mu_adt\\
&=-\iint_{P_1^+(0,1)} (|U|^2-1) {\rm{div}}_X(y^a\langle U, DU\rangle\eta^2)
+\iint_{P_1^+(0,1)} D(|U|^2-1)\langle U, U-(U)^a_{(0,1), 2r, +}\rangle D\eta^2\mu_adt\\
&=-\iint_{P_1^+(0,1)} (|U|^2-1) \Big(\langle U, \partial_t U\rangle \eta^2+ |DU|^2\eta^2+\langle U, DU\rangle D\eta^2\Big)\mu_adt\\
&\quad+\iint_{P_1^+(0,1)} D(|U|^2-1)\langle U, U-(U)^a_{(0,1), 2r, +}\rangle D\eta^2\mu_adt\\
&=IV_1+IV_2.
\end{align*} 
Applying \eqref{dist-est1} and H\"older's inequality, we can bound $IV_1$  by
\begin{align*}
\big|IV_1\big|&\le \big\||U|^2-1\big\|_{L^\infty(P^+_\frac14(0,1))}\iint_{P^+_1(0,1)}|DU|^2\eta^2\mu_adt\\
&+2\iint_{P^+_1(0,1)} ||U|^2-1|(|\partial_t U|\eta^2+|DU||D\eta^2|)\mu_adt\\
&\le \frac{\varepsilon_0}2 \iint_{P^+_{2r}(0,1)}(r^2|\partial_t U|^2+|DU|^2)\mu_adt+C\varepsilon_0^{-1} r^{-2} \iint_{P^+_{2r}(0,1)}(|U|^2-1)^2\mu_adt,
\end{align*}
while $IV_2$ can be bounded by
\begin{align*}
\big|IV_2\big|\le \frac{\varepsilon_0}2 \iint_{P_{2r}^+(0,1)}|DU|^2\mu_adt+C\varepsilon_0^{-1}r^{-2}\iint_{P_{2r}^+(0,1)}|U-(U)^a_{(0,1), 2r, +}|^2\mu_adt.
\end{align*}
Substituting the estimates of $IV_1$ and $IV_2$ into that of $IV$, we arrive at 
\begin{align*}
&\big|IV\big|\le \varepsilon_0 \iint_{P^+_{2r}(0,1)}\big(r^2|\partial_t U|^2+|DU|^2\big)\mu_adt
+C\varepsilon_0^{-1} r^{-2} \iint_{P^+_{2r}(0,1)}|U-(U)^a_{(0,1), 2r, +}|^2\mu_adt\\
&\ \ \qquad +C\varepsilon_0^{-1}r^{-2} \iint_{P^+_{2r}0,1)}(|U|^2-1)^2\mu_adt.
\end{align*}

The term $V$ can be handled similarly to $IV$. Integrating by parts,  we obtain
\begin{align*}
&|V|=\\
&\Big|\iint_{P_1^+(0,1)} (|U|^2-1)\Big(\partial_t U \cdot (U-(U)^a_{(0,1), 2r, +}) \eta^2+|DU|^2\eta^2+DU(U-(U)^a_{(0,1), 2r, +})\cdot D\eta^2\Big)\mu_adt\Big|\\
&\le  \big\||U|^2-1\big\|_{L^\infty(P^+_\frac14(0,1))}\iint_{P_{2r}^+(0,1)} \Big((r|\partial_t U|+|DU|)r^{-1}|U-(U)^a_{(0,1), 2r, +}|+|DU|^2\Big)\mu_adt\\
&\le \varepsilon_0  \iint_{P_{2r}^+(0,1)}\big(r^2|\partial_t U|^2+|DU|^2\big)\mu_adt+C\varepsilon_0 r^{-2} \iint_{P^+_{2r}(0,1)}|U-(U)^a_{(0,1), 2r, +}|^2\mu_adt.
\end{align*}

In order to estimate $VI$, we first apply \eqref{alg1} of Lemma \ref{alg0} to write
\begin{align*}
VI&=\int_{1-4r^2}^{1+4r^2}\int_{B_{2r}^+(0)}DU^i\cdot\Omega^{ki}\big(U^k-(U^k)^a_{(0,1), 2r, +}\big)\eta^2\mu_adt\\
&=-\int_{1-4r^2}^{1+4r^2}\int_{B_{2r}^+(0)}(U^i-(U^i)^a_{(0,1), 2r, +})\big(U^i\partial_tU^k-U^k\partial_t U^i\big))\big(U^k-(U^k)^a_{(0,1), 2r, +}\big)\eta^2\mu_adt\\
&\ \ \ -\int_{1-4r^2}^{1+4r^2}\int_{B_{2r}^+(0)}(U^i -(U^i)^a_{(0,1), 2r, +})\Omega^{ki}D\big(\big(U^k-(U^k)^a_{(0,1), 2r, +}\big)\eta^2\big)\mu_adt\\
&=VI_1+VI_2.
\end{align*}

By H\"older's inequality, $IV_1$ can be bounded by
\begin{align*}
\big|VI_1\big|\le \frac{\varepsilon_0}2\iint_{P_{2r}^+(0,1)} r^2|\partial_t U|^2\mu_adt
+C\varepsilon_0^{-1} r^{-2} \iint_{P^+_{2r}(0,1)}|U-(U)^a_{(0,1), 2r, +}|^4\mu_adt.
\end{align*}

The estimate of main term $VI_2$ requires Hardy and BMO duality. 
For a.e. $t\in (1-4r^2, 1+4r^2)$, observe that $y^a\Omega^{ki}(\cdot, t)\in L^2(B_{2r}^+(0),\mu_{-a})$ with weight $\mu_{-a}=y^{-a}dX$.
Since $y^{-a}$ is a $A_2$-weight,  the Hodge decomposition guarantees that 
there exist 
$$G^{ki}(t)\in H^1(B_{2r}^+(0), \mu_{-a})\  \ {\rm{and}}\  \ H^{ki}(t)\in H^1_0\big((B_{2r}^+(0), \mu_{-a});
\wedge^2\mathbb R^{m+1}\big)$$ 
such that
\begin{equation}\label{hodge1}
y^a\Omega^{ki}(\cdot, t) =D G^{ki}(t)+{\rm{curl}} H^{ki}(t) \ \ {\rm{in}}\ \  B_{2r}^+(0),
\end{equation}
and 
\begin{align}\label{h1-est}
&\int_{B_{2r}^+(0)} |DG^{ki}(t)|^2\mu_{-a} +\int_{B_{2r}^+(0)} |DH^{ki}(t)|^2\mu_{-a} \nonumber\\
&\le C\int_{B_{2r}^+(0)} (y^a|\Omega^{ki}(\cdot, t)|)^2(X, t)\mu_{-a}=C\int_{B_{2r}^+(0)} |\Omega^{ki}(\cdot, t)|^2(X, t)\mu_a.
\end{align} 
In view of \eqref{alg1}, $G^{ki}(t)$ solves the Neumann boundary value problem
\begin{equation}\label{g-equ}
\begin{cases}
-\Delta_X G^{ki}(t)= y^a\big(U^i\partial_t U^k-U^k\partial_t U^i\big)(\cdot, t) & \ \ \ {\rm{in}}\ \ B_{2r}^+(0),\\
\displaystyle\qquad\frac{\partial G^{ki}(t)}{\partial y}=0 & \ \ \ {\rm{on}}\ \ \partial B_{2r}^+(0)\cap \partial \mathbb{R}^{m+1}_+.
\end{cases}
\end{equation}  
By the $y^{-a}$-weighted $W^{2,2}$-estimate, we have that $G^{ki}\in H^2(B_{\frac{3r}2}^+(0), \mu_{-a})$, and
\begin{align}\label{h2-est}
&\int_{B_{\frac{3r}2}^+(0)}\Big(|DG^{ki}|^2(X,t)+r^2|D^2G^{ki}|^2(X, t)\Big)\mu_{-a}\nonumber\\
&\ \le C\int_{B_{2r}^+(0)}\Big(|DG^{ki}|^2(X,t)+r^2y^{2a}|\partial_t U|^2(X, t)\Big)\mu_{-a}\nonumber\\
&\ \le C\int_{B_{2r}^+(0)}\Big(|y|^{2a}|\Omega^{ki}|^2(X,t)+r^2y^{2a}|\partial_t U|^2(X, t)\Big)\mu_{-a}\nonumber\\
&\ \le C\int_{B_{2r}^+(0)}\Big(|DU|^2(X,t)+r^2|\partial_t U|^2(X, t)\Big)\mu_a.
\end{align} 
By the $y^a$-weighted Sobolev inequality \eqref{weight-sobolev}, \eqref{h2-est} yields
\begin{align*}
\big\|DG^{ki}(t)\big\|_{L^{2_a^*}(B_{\frac{3r}2}^+(0,1),\  \mu_{-a})}&\le C
\Big(\|D^2G^{ki}(t)\|_{L^2(B_\frac{3r}2^+(0),\ \mu_{-a})}+r^{-1}\|DG^{ki}(t)\|_{L^2(B_\frac{3r}2^+(0),\ \mu_{-a})}\Big)\\
&\le Cr^{-1} \Big(\|DU(\cdot, t)\|_{L^2(B_{2r}^+(0),\ \mu_a)}+r\|\partial_t U(\cdot, t)\|_{L^2(B_{2r}^+(0),\ \mu_a)}\Big),
\end{align*}
where 
$$
2_a^*=\begin{cases}
\displaystyle\frac{2(m+1-a)}{m-1-a}, & \ m\ge 2 \ {\rm{or}}\  m=1 \ {\rm{and}} \ a<0,\\
 \infty,  & \ m=1 \ {\rm{and}} \ a>0.
 \end{cases}
 $$

Combining condition \eqref{small1}, the interior H\"older continuity of $U$, the energy monotonicity inequality \eqref{local_monotonicity1}, and the local energy inequality
\eqref{energy_ineq}
 and the $y^a$-weighted Poincar\'e inequality, there exists $\theta_0\in (0,1)$, depending only on $\varepsilon_0, t_0,E_0$,
 such that 
\begin{align}\label{weighted-sobolev}
\frac{1}{\displaystyle |P_r^+(z)|_a} \iint_{P_r^+(z)}|U-(U)^a_{z, r,+}|^2\mu_adt
&\le \displaystyle\frac{Cr^2}{\displaystyle |P_r^+(z)|_a} \iint_{P_r^+(z)}(r^2|\partial_t U|^2+|DU|^2)\mu_adt\nonumber\\
&\le C\iint_{P_{\frac12}^+(0,1)}(|\partial_t U|^2+|DU|^2)\mu_adt\nonumber\\
&\le C\varepsilon_0^2
\end{align}
holds for all $z\in P_{\frac12}^+(0,1)$ and $0<r<\theta_0$. This implies $U\in {\rm{BMO}}(P_{\theta_0}^+(0,1), \mu_adt)$ with bound
\begin{align}\label{weighted-BMO}
\big[u\big]_{{\rm{BMO}}(P_{\theta_0}^+(0,1),\ \mu_adt)}\le C\varepsilon_0.
\end{align}
Since $y^a$ is a $A_2$-weight, we can apply the weighted form of John-Nirenberg's inequality \eqref{weight-JN0} to obtain that for any $z\in P_{\theta_0}^+(0,1)$ and $0<r<{\theta_0}$,
$\big(U-(U)^a_{z, r, +}\big)\in L^p\big(P_{r}^+(z), \mu_adt\big)$ for all $p\in (1,\infty)$, satisfying
\begin{align}\label{weighted-JN}
\displaystyle\frac{1}{\displaystyle|P_r^+(z)|_a} \iint_{P_r^+(z)}|U-(U)^a_{z, r,+}|^p\mu_adt
\le C(m,p) \big[u\big]^p_{{\rm{BMO}}(P_{\theta_0}^+(0,1),\ \mu_adt)}
\le C(m, p)\varepsilon_0^p.
\end{align}
Decomposing $IV_2$ by \eqref{hodge1}, we have 
\begin{align*}
-VI_2&=\int_{1-4r^2}^{1+4r^2}\int_{B_{2r}^+(0)}(U^i -(U^i)^a_{(0,1), 2r, +})\Omega^{ki}D\big(\big(U^k-(U^k)^a_{(0,1), 2r, +}\big)\eta^2\big)\mu_adt\\
&=\int_{1-4r^2}^{1+4r^2}\int_{B_{2r}^+(0)} (U^i -(U^i)^a_{(0,1), 2r, +})DG^{ki}D\big(\big(U^k-(U^k)^a_{(0,1), 2r, +}\big)\eta^2\big)\\
&\quad+\int_{1-4r^2}^{1+4r^2}\int_{B_{2r}^+(0)} (U^i -(U^i)^a_{(0,1), 2r, +}){\rm{curl}} H^{ki}D\big(\big(U^k-(U^k)^a_{(0,1), 2r, +}\big)\eta^2\big)\\
&=VI_{21}+VI_{22}.
\end{align*}
Using H\"older's inequality with conjugate exponent $(2_a^{*})'=\frac{2_a^*}{2_a^*-1}$,  $VI_{21}$ is  bounded by
\begin{align*}
\big|VI_{21}\big|
&\le  \Big(\int_{1-4r^2}^{1+4r^2} \big\|DG\big\|_{L^{2_a^*}(B_{2r}^+(0), \mu_{-a})}^2\,dt\Big)^\frac12\\
&\quad\cdot\Big(\int_{1-4r^2}^{1+4r^2}\big\|D\big((U-(U)^a_{(0,1), 2r, +})\eta^2\big)\big(U-(U)^a_{(0,1), 2r, +}\big)\big\|_{L^{(2_a^*)'}(B_{2r}^+(0), \mu_a)}^2\Big)^\frac12\\
&\le Cr^{-1}\Big(\iint_{P_{2r}^+(0,1)} (r^2|\partial_t U|^2+|DU|^2)\mu_adt\Big)^\frac12\\
&\quad\cdot \Big(\int_{1-4r^2}^{1+4r^2} \Big\|r^{-1}|U-(U)^a_{(0,1), 2r, +}|^2+|DU||U-(U)^a_{(0,1), 2r,+}|\Big\|_{L^{(2_a^*)'}(B_{2r}^+(0), \mu_a)}^2\,dt\Big)^\frac12\\
&\le Cr^{-1}\Big(\iint_{P_{2r}^+(0,1)} (r^2|\partial_t U|^2+|DU|^2)\mu_adt\Big)^\frac12\\
&\quad\cdot \Big(\int_{1-4r^2}^{1+4r^2} \Big\|r^{-1}|U-(U)^a_{(0,1), 2r, +}|^2+|DU||U-(U)^a_{(0,1), 2r,+}|\Big\|_{L^{(2_a^*)'}(B_{2r}^+(0), \mu_a)}^2\,dt\Big)^\frac12\\
&\le C\varepsilon_0 \iint_{P_{2r}^+(0,1)}(r^2|\partial_t U|^2+|DU|^2)\mu_adt 
+C\varepsilon_0^{-1}r^{-2}\iint_{P_{2r}^+(0,1)} |U-(U)^a_{(0,1), 2r, +}|^4\mu_adt.
\end{align*}

Finally, we estimate $VI_{22}$ using the Hardy--BMO duality. Extending $H^{ki}(t)$ by $0$ outside $B_{2r}^+(0)$
and $U$  evenly across $\partial\R^{m+1}_+$, we apply the $y^a$-weighted parabolic div-curl Lemma \ref{hardy} from Appendix to obtain
\begin{align*}
\big|VI_{22}\big|&= \Big|\int_{1-4r^2}^{1+4r^2}\int_{\R^{m+1}}(U^i-(U^i)^a_{(0,1), 2r,+})\cdot ({\rm{curl}} H^{ki})\cdot D\big(\big(U^k-(U^k)^a_{(0,1), 2r}\big)\eta^2\big)\Big|\\
&\le C\big\|D\big(\big(U^k-(U^k)^a_{(0,1), 2r}\big)\eta^2\big) \cdot {\rm{curl}}H^{ki}\big\|_{\mathcal{H}^1(\R^{m+1}\times \R)}\big[U^i\big]_{{\rm{BMO}}\big(P_{2r}^+(0,1)\big)}\\
&\le C\big\|DU^i\big\|_{L^2\big(P_{2r}^+(0,1),\ \mu_adt\big)}\big\|DH^{ki}\big\|_{L^2\big(P_{2r}^+(0,1),\ \mu_{-a}dt\big)} 
\big[U^i\big]_{{\rm{BMO}}\big(P_{2r}^+(0,1)\big)}\\
&\le C\varepsilon_0\iint_{P_{2r}^+(0,1)}(r^2|\partial_t U|^2+|DU|^2)\mu_adt
\end{align*}
where we have used \eqref{h1-est} and the BMO estimate (see \eqref{bmo-bmoa})
\begin{align*}
\Big[U\Big]_{{\rm{BMO}}\big(P_{2r}^+(0,1)\big)}
&\le \Big[U\Big]_{{\rm{BMO}}\big(P_{\theta_0}^+(0,1)\big)}\le C\Big[U\Big]_{{\rm{BMO}}\big(P_{\theta_0}^+(0,1),\ \mu_adt\big)}\\
&\le C\Big(\big\|DU\big\|_{M^{2,2}(P_\frac12^+(0,1),\ \mu_adt)}+\big\|\partial_tU\big\|_{M^{2,4}(P_\frac12^+(0,1),\ \mu_adt)}\Big)\\
&\le C\varepsilon_0.
\end{align*}
Combing the bounds for $I$ through $VI$ yields \eqref{hybrid2}. This completes the proof
\qed
 
\section{Energy monotonicity inequality and local energy inequality}

In this section, motivated by the definition of stationarity of heat flow of harmonic maps, we introduce the notation
of stationary heat flow of $s$-harmonic maps and derive their fundamental properties.

For convenience, we first recall the Definition 1.2. Let $\xi(X,t)\in C^\infty_0\big(\overline{\mathbb R^{m+1}_+}\times (0,\infty), \mathbb R^{m+1}\big)$
with $\xi^{m+1}\ge 0$ and $\theta (X,t)\in C^\infty_0\big(\overline{\mathbb R^{m+1}_+}\times (0,\infty)\big)$ with $\theta\ge 0$, and let
$U\in L^\infty\big(\R_+;  \dot{H}^1((\R^{m+1}_+,\mu_a), \R^{\ell})\big)\cap \dot{H}^1\big(\R_+;  L^2((\R^{m+1}_+,\mu_a), \R^{\ell})\big).
$
For a sufficiently small $\tau>0$, define the inner variation $U^\tau$ of $U$ by
$$
U^\tau(X,t)=U\big(X+\tau\xi(X, t), t+\tau\theta(X,t)\big), \ (X,t)\in\overline{\mathbb R^{m+1}_+}\times [0,\infty).
$$

\begin{definition} \label{stationary10} A weak solution $U:{\mathbb R^{m+1}_+}\times (0,\infty)\to \mathbb R^{\ell}$ of \eqref{eq:fhf3}
is called a stationary heat flow of $s$-harmonic maps, if the following inequality
\begin{equation}\label{stationary1}
2\iint_{\mathbb R^{m+1}_+\times \mathbb R_+} \partial_tU \cdot\frac{d}{d\tau}\big|_{\tau=0} U^\tau\mu_adt
+\frac{d}{d\tau}\big|_{\tau=0}\iint_{\mathbb R^{m+1}_+\times\mathbb R_+} \big|DU^\tau\big|^2\mu_adt\le 0.
\end{equation}
holds for all $\xi(X,t)\in C^\infty_0\big(\overline{\mathbb R^{m+1}_+}\times (0,\infty), \mathbb R^{m+1}\big)$,
with $\xi^{m+1}\ge 0$ and $\xi^{m+1}((x,0),t)=0$ for $(x,t)\in\mathbb R^m\times (0,\infty)$,
and $0\le \theta (X,t)\in C^\infty_0\big(\overline{\mathbb R^{m+1}_+}\times (0,\infty)\big)$.
\end{definition} 

\medskip
Since $\xi(\cdot, t): \partial\mathbb R^{m+1}_+\to \partial\mathbb R^{m+1}_+$ for $t>0$, it is clear
that $U^\tau(x,0,t)\in \mathbb S^{\ell-1}$ for $(x,t)\in\mathbb R^m\times (0,\infty)$, rendering $U^\tau$ an admissible
inner variation of $U$.  A direct calculation yields
$$
\frac{d}{d\tau}\big|_{\tau=0} U^\tau=\xi\cdot DU +\theta \partial_tU,
$$
so that condition \eqref{stationary1} takes the form
\begin{align}\label{stationary2}
\iint_{\mathbb R^{m+1}_+\times \mathbb R_+} \partial_tU\cdot\big(\xi\cdot DU +\theta \partial_tU\big)\mu_adt
+\iint_{\mathbb R^{m+1}_+\times\mathbb R_+} DU\cdot D\big(\xi\cdot DU +\theta \partial_tU\big)\mu_adt\le 0.
\end{align}

For any $0<t_1<t_2<\infty$, if we set $\theta\equiv 0$ and choose $\xi\in C_0^\infty(\overline{\R^{m+1}_+}\times (t_1,t_2))$
with $\xi^{m+1}=0$ on $\partial\R^{m+1}_+\times (0,\infty)$ (and testing with both $\xi$ and  $-\xi$) in \eqref{stationary2}, 
and integration by parts together with the boundary condition of $U$ on $\partial\R^{m+1}_+\times (0,\infty)$ yields the 
the stationarity identity 
\begin{align}\label{stationary3}
\iint_{\mathbb R^{m+1}_+\times [t_1, t_2]} y^a\partial_tU\cdot\big(\xi\cdot DU\big)
-\frac12 {\rm{div}}_X(y^a\xi)|DU|^2
+y^aDU\otimes DU: D\xi=0.
\end{align}
On the other hand, setting $\xi\equiv 0$ and choosing $0\le \theta\in C_0^\infty(\R^{m+1}_+\times [t_1,t_2])$  in \eqref{stationary2}
results
\begin{align}\label{stationary4}
&\iint_{\mathbb R^{m+1}_+\times \{t_2\}} |DU|^2\theta\mu_a
-\iint_{\mathbb R^{m+1}_+\times \{t_1\}} |DU|^2\theta\mu_a\nonumber\\
&\le \iint_{\mathbb R^{m+1}_+\times [t_1, t_2]} \big(|DU|^2\partial_t\theta-2D\theta\cdot DU\cdot \partial_tU-2|\partial_t U|^2\theta\big)\mu_adt.
\end{align}

\medskip
For $X_0=(X_0,0)\in\partial\mathbb R^{m+1}_+$ and $t_0>0$, let $G^a_{(X_0, t_0)}$  be the backward kernel  of the weighted parabolic operator
$\mathcal{L}_a V:=|y|^a \partial_t V-{\rm{div}}(|y|^a DV)$ on $\R^{m+1}$,  defined by
\begin{equation}
G^a_{(X_0, t_0)}(X,t)=\frac{1}{\Gamma(\frac{1-a}2)(4\pi)^{\frac{m}2} |t-t_0|^{\frac{m+1+a}2}}\exp(-\frac{|X-X_0|^2}{4|t-t_0|}), \  t<t_0.
\end{equation}
In particular, when $s=1/2$,  $a=0$ and $G^0_{(X_0,t_0)}$ becomes the backward heat kernel $G_{(X_0,t_0)}$ on $\R^{m+1}$.
For $0<R\le \frac{\sqrt{t_0}}2$, set 
\begin{align*}T_R^+(t_0)&=\big\{ (X,t)\in\R^{m+1}_+\times\R_+:  t_0-4R^2< t< t_0-R^2 \big\}, \\
\partial^+ T_R^+(t_0)&=\big\{ (x, t)\in\R^m\times\R_+:t_0 -4R^2<t<t_0-R^2 \big\} ,\end{align*}
and abbreviate $T_1^+=T_1^+(0)$ and $\partial^+ T_1^+=\partial T_1^+(0)$.

\begin{lemma}\label{lem-mono} 
If $U$ is a stationary weak solution of heat flow of $s$-harmonic maps and $a=1-2s$
 then, for every $Z_0=(X_0,t_0)$ with $X_0\in\partial\R^{m+1}_+$ and $t_0>0$, the renormalized energies
\begin{align*}
\mathcal{D}^a(U, Z_0, R)&:=\frac{R^2}2\int_{\R^{m+1}_+\times \{t_0-R^2\}}|DU|^2G^a_{(X_0,t_0)}(X,t)\mu_a
\end{align*}
and
\begin{align*}
\mathcal{E}^a(U, Z_0, R)&:=  \frac12\iint_{T_R^+(t_0)}|DU|^2G^a_{(X_0,t_0)}(X,t)\,\mu_adt
\end{align*}
are monotone nondecreasing with respect to $R$, satisfying
\begin{align}\label{mono1}
\mathcal{D}^a(U, Z_0, r)\le \mathcal{D}^a(U, Z_0, R), \ 0<r\le R<\sqrt{t_0},
\end{align}  
and
\begin{align}\label{mono2}
\mathcal{E}^a(U, Z_0, r)\leq \mathcal{E}^a(U, Z_0, R), \ \ 0<r\le R< \frac12\sqrt{t_0}.
\end{align}
\end{lemma}
\pf Without loss of generality, assume $X_0=(0,0)$ and $t_0=1$, and write $G^a(X,t)=G^a_{(X_0, t_0)}(X,t)$. 
Note that testing with $-\xi$ makes \eqref{stationary3} an equality.  Select $\xi=XG^a(X,t)$. Then $\xi^{m+1}$ is nonnegative and satisfies
$\xi^{m+1}=0$ on $\partial\mathbb R^{m+1}_+\times (0,\infty)$ \footnote{Strictly speaking, we need to multiply $\xi$ by a sequence of cut-off functions 
$\eta_n\in C_0^\infty(\overline{\R^{m+1}_+})$, with $\eta_n\rightarrow 1$. But, it is not hard to check that the error terms resulting from $\eta_n$ tend to zero as $n\to\infty$. Thus we simply
choose $\eta_n\equiv 1$.}. 
Direct computations yield
\begin{align*}
{\rm{div}}_X(y^a\xi)&=y^a\big(m+1+a -\frac{|X|^2}{2(1-t)}\big)G^a, \\ 
\nabla_{X_i}\xi^j &=\begin{cases} \big(\delta_{ij}-\frac{X_i X_j}{2(1-t)}\big) G^a, & \ 1\le i, j\le m,\\
-\frac{yX_i}{2(1-t)} G^a, & \ 1\le i\le m, \ j=m+1,\\
-\frac{yX_j}{2(1-t)}G^a, &\  i=m+1, \ 1\le j\le m,\\
  \big(1-\frac{y^2}{2(1-t)}\big)G^a,  &\ i=j=m+1.
\end{cases}
\end{align*}
For $0<r<R<1$, setting $t_1=1-R^2$ and $t_2=1-r^2$ and substituting $\xi$ into \eqref{stationary3} yields
\begin{align}\label{mono3}
0&=\iint_{\mathbb R^{m+1}_+\times [1-R^2, 1-r^2]} \partial_tU\cdot(X\cdot DU)G^a\mu_adt\nonumber\\
& \ \ +\iint_{\mathbb R^{m+1}_+\times [1-R^2, 1-r^2]} \Big[\big(-\frac{m-1+a}2 +\frac{|X|^2}{4(1-t)}\big)|DU|^2
-\frac{|X\cdot DU|^2}{2(1-t)}\Big]G^a\mu_adt.
\end{align}
Next, choosing $\theta=(1-t)G^a$, we find $\theta_t=\frac{m-1+a}2G^a-\frac{|X|^2}{4(1-t)} G^a$. Substituting $\theta$ into
\eqref{stationary4} yields
\begin{align}\label{mono4}
&\iint_{\mathbb R^{m+1}_+\times [1-R^2, 1-r^2]} \big(-2(1-t)|\partial_tU|^2
+X\cdot DU\cdot \partial_tU+(\frac{m-1+a}2-\frac{|X|^2}{4(1-t)})|DU|^2\big)G^a\mu_adt\nonumber\\
&\ge \int_{\mathbb R^{m+1}_+\times \{1-r^2\}} |DU|^2 r^2G^a\mu_a
-\int_{\mathbb R^{m+1}_+\times \{1-R^2\}} |DU|^2 R^2G^a\mu_a.
\end{align}
Adding \eqref{mono3} and \eqref{mono4}, we obtain 
\begin{align}\label{mono5}
&R^2\int_{\mathbb R^{m+1}_+\times \{1-R^2\}} |DU|^2 G^a\mu_a
-r^2\int_{\mathbb R^{m+1}_+\times \{1-r^2\}}|DU|^2G^a\mu_a\nonumber\\
&\ge \iint_{\mathbb R^{m+1}_+\times [1-R^2, 1-r^2]} \big(2(1-t)|\partial_tU|^2
-2X\cdot DU\cdot \partial_tU+\frac{1}{2(1-t)}|X\cdot DU|^2\big)G^a\mu_adt\nonumber\\
&=\iint_{\mathbb R^{m+1}_+\times [1-R^2, 1-r^2]}\frac{|X\cdot DU+2(t-1) \partial_tU|^2}{2(1-t)}G^a\mu_adt.
\end{align}
This proves \eqref{mono1}. Moreover, observing that
\[
\mathcal{E}^a(U, Z_0, R)=2\int_R^{2R}\frac{\mathcal{D}^a(U, Z_0, \tau)}{\tau}\,d\tau,
\]
so that differentiating with respect to $R$ implies
\[
\frac{d}{dR}\mathcal{E}^a(U, Z_0, R)=\frac{2}{R} \big(D^a(U, Z_0, 2R)-D^a(U, Z_0, R)\big)\ge 0.
\] 
This implies \eqref{mono2}. 
\qed

Next,  we obtain a local version of $y^a$-weighted energy monotonicity inequality.

\begin{lemma}\label{local_monotonicity} Let $U$ be a stationary solution of heat flow of $s$-harmonic maps in $\R^{m+1}_+$ and $a=1-2s$.
For any $\varepsilon>0$, there exists $K(\varepsilon)\sim \sqrt{|\ln \varepsilon|}$ such that for $Z_0=(X_0,t_0)\in\partial\R^{m+1}_+\times (0,\infty)$ 
and $0<R<\frac{\sqrt{t_0}}{2K(\varepsilon)}$,  
\begin{equation}\label{local_monotonicity1} 
\rho^{-(m+1+a)}\iint_{P_{\rho}^+(Z_0)} |DU|^2\mu_adt\le 2\ln(\frac52)\varepsilon t_0^{\frac{1-m-a}2}E_0
+R^{-(m+1+a)}\iint_{P_{K(\varepsilon) R}^+(Z_0)}|DU|^2\mu_adt.
\end{equation}
holds for $0<\rho<\frac{R}2$,  where $E_0=\displaystyle\int_{\R^{m+1}_+} y^a|DU|^2(\cdot,0)$.
\end{lemma}
\pf  
By the definition of $G^a$,  there is $c=c(m,a)>0$ such that
\[
\rho^{-(m+1+a)}\le cG^a_{(X_0, t_0+2\rho^2)}(X,t), \  \forall (X,t)\in P_\rho^+(Z_0).
\]
Since $P_\rho^+(Z_0)\subset T^+_\rho(t_0+2\rho^2)$, by \eqref{mono2} we have that
\begin{align}\label{local_monotonicity1.1}
\rho^{-(m+1+a)}\iint_{P_{\rho}^+(Z_0)} |DU|^2\mu_adt&\le c \iint_{P_{\rho}^+(Z_0)} |DU|^2G^a_{(X_0, t_0+2\rho^2)}\mu_adt\nonumber\\
&\le c\iint_{T^+_\rho(t_0+2\rho^2)} |DU|^2G^a_{(X_0, t_0+2\rho^2)}\mu_adt\nonumber\\
&\le  c\iint_{T^+_R(t_0+2\rho^2)} |DU|^2G^a_{(X_0, t_0+2\rho^2)}\mu_adt.
\end{align}
Now, we claim that 
\begin{align}\label{local_monotonicity2}
\iint_{T^+_R(t_0+2\rho^2)} |DU|^2G^a_{(X_0, t_0+2\rho^2)}\mu_adt
&\le 2\ln(\frac52) \varepsilon t_0^{\frac{1-m-a}2}E_0\nonumber\\
&\quad+R^{-(m+1+a)}\iint_{P_{K(\varepsilon) R}^+(Z_0)} |DU|^2\mu_adt.
\end{align} 
The proof  of \eqref{local_monotonicity2} is based on suitable modifications of \cite[Theorem 3]{Cheng} (see also \cite[Lemma 3.1] {Feldman}). 
Here we sketch it. 
First, as in the proof of \cite[Lemma 1]{Cheng}, if we choose $K(\varepsilon)\sim \sqrt{|\log\varepsilon|}$, 
then there holds, 
\begin{align*}
G^a_{(X_0, t_0+2\rho^2)}(X,t)\le\begin{cases} R^{-(m+1+a)} & \ (X,t)\in T_R^+(t_0+2\rho^2),\\
\varepsilon G^a_{(X_0, t_0+2\rho^2+R^2)}(X,t) & (X,t)\in T_R^+(t_0+2\rho^2) \ {\rm{with}}\ |X-X_0|\ge K(\varepsilon) R.
\end{cases}
\end{align*}
Indeed, for $(X,t)\in T_R^+(t_0+2\rho^2)$ with $|X-X_0|\ge K(\varepsilon) R$,
\begin{align*}
\frac{G^a_{(X_0, t_0+2\rho^2)}(X,t)}{G_{(X_0,t_0+2\rho^2+R^2)}^a(X,t)}
&= \big(\frac{t_0+2\rho^2+R^2-t}{t_0+2\rho^2-t}\big)^{\frac{m+1+a}2}\exp\Big(-\frac{|X-X_0|^2}{4(t_0+2\rho^2-t)}+
\frac{|X-X_0|^2}{4(t_0+2\rho^2+R^2-t)}\Big)\\
&\le 5^{\frac{m+1+a}{2}} \exp\Big(-\frac{R^2|X-X_0|^2}{4(t_0+2\rho^2-t)(t_0+2\rho^2+R^2-t)}\Big)\\
&\le 5^{\frac{m+1+a}{2}} \exp\Big(-\frac{K^2(\varepsilon)}{8}\Big)\leq \varepsilon.
\end{align*}
Decomposing the integration domain into $B_{K(\varepsilon)R}^+(X_0)$ and $T_R^+(t_0+2\rho^2)\setminus B_{K(\varepsilon)R}^+(X_0)$, we obtain 
\begin{align}\label{mono10}
&\iint_{T^+_R(t_0+2\rho^2)} |DU|^2G^a_{(X_0, t_0+2\rho^2)}\mu_adt\nonumber\\
&\le R^{-(m+1+a)} \int_{t_0+2\rho^2-4R^2}^{t_0+2\rho^2-R^2}\int_{B^+_{K(\varepsilon) R}(X_0)}|DU|^2\mu_adt
+\varepsilon \iint_{T_R^+(t_0+2\rho^2)}|DU|^2 G_{(X_0,t_0+2\rho^2+R^2)}^a\mu_adt\nonumber\\
&\le R^{-(m+1+1)} \iint_{P^+_{K(\varepsilon)R}(Z_0)} |DU|^2\mu_adt
+\varepsilon \iint_{T_R^+(t_0+2\rho^2)}|DU|^2 G_{(X_0,t_0+2\rho^2+R^2)}^a\mu_adt.
\end{align} 
As in \cite[page 172]{Cheng}, we can apply Lemma \ref{lem-mono} to estimate
\begin{align}\label{mono20}
&\iint_{T_R^+(t_0+2\rho^2)}|DU|^2 G_{(X_0,t_0+2\rho^2+R^2)}^a\mu_adt\nonumber\\
&=2\int_{t_0+2\rho^2-4R^2}^{t_0+2\rho^2-R^2}(t_0+2\rho^2+R^2-t)^{-1} \mathcal{D}^a\Big(U,  (X_0, t_0+2\rho^2+R^2), \sqrt{t_0+2\rho^2+R^2-t}\Big)\,dt\nonumber\\
&\le2\mathcal{D}^a\Big(U,  (X_0, t_0+2\rho^2+R^2), \sqrt{t_0+2\rho^2+R^2}\Big)\int_{t_0+2\rho^2-4R^2}^{t_0+2\rho^2-R^2}(t_0+2\rho^2+R^2-t)^{-1}\,dt\nonumber\\
&\le 2\ln(\frac52) (t_0+2\rho^2+R^2)^{\frac{1-m-a}2} E_0\le 2\ln(\frac52) t_0^{\frac{1-m-a}2} E_0.
\end{align} 
Combining \eqref{mono10} with \eqref{mono20} yields \eqref{local_monotonicity2}. Substituting \eqref{local_monotonicity2}
into \eqref{local_monotonicity1.1} completes the proof of \eqref{local_monotonicity1}.
\qed

\bigskip
The stationarity condition \eqref{stationary1} also leads to the local energy inequality. More precisely,
let $0\le\eta\in C^\infty_0(\mathbb R^{m+1}\times (0,\infty))$ be even w.r.t. $y$, i.e. $\eta((x, -y), t)=\eta(x, y, t)$
for $X=(x,y)\in\mathbb R^{m+1}$ and $t>0$. Then $\xi(X,t)=D\eta\in C^\infty(\overline{\mathbb R^{m+1}_+})$ 
satisfies $\xi^{m+1}((x,0), t) =0$ for $(x,t)\in\mathbb R^m\times (0,\infty)$. Substituting $\xi$ into \eqref{stationary3} and
replacing $\theta$ by $\eta$ in \eqref{stationary4} and summing the two resulting inequalities yields
\begin{align}\label{energy_ineq}
& \int_{\mathbb R^{m+1}_+\times \{t_2\}} |DU|^2\eta\mu_a
-\int_{\mathbb R^{m+1}_+\times \{t_1\}} |DU|^2\eta\mu_a+2\iint_{\mathbb R^{m+1}_+\times [t_1, t_2]} |\partial_tU|^2\eta\mu_adt
\nonumber\\
&\le \iint_{\mathbb R^{m+1}_+\times [t_1, t_2]} \Big[|DU|^2\big(y^a\partial_t\eta+{\rm{div}}_X(y^a D\eta)\big)
-2y^aDU\otimes DU: D^2\eta \Big]dXdt.
\end{align} 

As an immediate consequence of \eqref{energy_ineq}, we establish the local $y^a$-energy inequality.
\begin{lemma}\label{local_energy_ineq}  If $U$ is a stationary solution of heat flow of $s$-harmonic maps in $\R^{m+1}_+$ and $a=1-2s$
 then 
 \begin{equation} \label{local_energy_ineq1}
 \iint_{P_R^+(Z_0)} |\partial_t U|^2\mu_adt\le C R^{-2}\int_{P_{2R}^+(Z_0)} |DU|^2\mu_adt
 \end{equation}
 holds for every $Z_0=(X_0,t_0)\in\partial{\R^{m+1}_+}\times (0,\infty)$ and $0<R<\frac{\sqrt{t_0}}4$.
\end{lemma}
\pf By Fubini's theorem there exists $t_*\in (t_0-R^2, t_0-4R^2)$ such that 
\[
\int_{B^+_{2R}(X_0)} |DU|^2(X, t_*)\mu_a\le \frac{1}{R^2}\iint_{P_{2R}^+(Z_0)}|DU|^2\mu_adt.
\]
Choose a parabolic cut-off function $\eta\in C^\infty_0(P_{2R}(Z_0))$ satisfying $\eta(x,-y,t)=\eta(x,y,t)$,
$0\le\eta\le 1$,  $\eta\equiv 1$ in $P_R(Z_0)$, and $|D^2\eta|+|D\eta|^2+|\partial_t \eta|\le \frac{24}{R^2}$.
Since $\frac{\partial \eta}{\partial y}\big|_{y=0}=0$, we have
\[
|\frac{\partial\eta}{\partial y}|\le |y|\big\|\frac{\partial^2 \eta}{\partial y^2}\big\|_{L^\infty}\le CR^{-2} |y|.
\]
With help of this, we can apply \eqref{local_energy_ineq1} with $t_=t_*$ and $t_2=t_0+R^2$ to get
\begin{align*}
\iint_{P^+_R(Z_0)}|\partial_tU|^2\mu_adt&\le \iint_{\R^{m+1}_+\times [t_*, t_0+R^2]} |\partial_tU|^2\eta\mu_adt\\
&\le \int_{\R^{m+1}_+\times \{t_*\}}|DU|^2\eta\mu_a\\
&+ \iint_{\mathbb R^{m+1}_+\times [t_*, t_0+R^2]} [|DU|^2\big(y^a\partial_t\eta+{\rm{div}}_X(y^a D\eta)\big)
-2y^aDU\otimes DU: D^2\eta]dXdt\\
&\le \int_{B_{2R}^+(X_0)\times \{t_*\}}|DU|^2\mu_a+C\iint_{P_{2R}^+(Z_0)}\big(R^{-2}DU|^2+y^{-1}|\frac{\partial \eta}{\partial y}||DU|^2\big)\mu_adt\\
&\le CR^{-2}\iint_{P_{2R}^+(Z_0)}|DU|^2\mu_adt.
\end{align*} 
This yields \eqref{local_energy_ineq1}. The proof is complete.
\qed

\bigskip

\section{Proof of Theorem \ref{partial_reg}: H\"older continuity}

In this section, we will prove the following partial regularity theorem for stationary solutions of the heat flow of $s$-harmonic maps. 
A key step to prove Theorem \ref{partial_reg} is  the following $\epsilon_0$-decay Lemma.

\begin{lemma}\label{decay1} There exist $\epsilon_0>0$ and $\theta_0\in (0,1)$ such that if $U$ is a stationary solution of heat flow of $s$-harmonic maps in $\R^{m+1}_+$
satisfying, for $(X_0,t_0)\in\partial\mathbb R^{m+1}_+\times (0,\infty)$ and $0<r<\frac{\sqrt{t_0}}{4}$,
\begin{equation}\label{decay2}
(2r)^{-(m+1+a)}\iint_{P_{2r}^+(X_0,t_0)} |DU|^2 \mu_adt\leq \epsilon_0^2,
\end{equation}
then
\begin{equation}\label{decay3}
(\theta_0 r)^{-(m+1+a)}\iint_{P_{\theta_0 r}^+(X_0,t_0)}|DU|^2 \mu_adt\le \frac12 (2r)^{-(m+1+a)} \iint_{P_{2r}^+(X_0,t_0)} |DU|^2\mu_adt
\end{equation}
\end{lemma}
\pf We may assume $X_0=0$ and $t_0=1$. We argue it by contradiction. Suppose the conclusions were false. Then  for any $\theta\in (0,1)$ there exist $r_i>0$ and
$\varepsilon_i\to 0$ such that
\begin{equation}\label{decay4}
(2r_i)^{-(m+1+a)}\iint_{P_{2r_i}^+(0,1)} |DU|^2\mu_adt=\varepsilon_i^2,
\end{equation}
then
\begin{equation}\label{decay5}
(\theta r_i)^{-(m+1+a)}\iint_{P_{\theta r_i}^+(0,1)}|DU|^2\mu_adt> \frac12 (2r_i)^{-(m+1+a)} \iint_{P_{2r_i}^+(0,1)} |DU|^2\mu_adt.
\end{equation}
Define the blow-up sequence $U_i(X,t)=\frac{U(r_i X, 1+r_i^2t)-U(0,1)}{\epsilon_i}:\overline{\mathbb R^{m+1}_+}\times (-4,4)\to\mathbb R^{\ell}$.
Then we have that   
\begin{equation}\label{decay6}
2^{-(m+1+a)}\iint_{P_{2}^+(0,0)} |DU_i|^2\mu_adt=1,
\end{equation}
but
\begin{equation}\label{decay7}
\theta^{-(m+1+a)}\iint_{P_{\theta}^+(0,0)}|DU_i|^2\mu_adt> \frac12.
\end{equation}
Since $U$ solves \eqref{eq:fhf3}, it follows that $U_i$ satisfes
\begin{equation}\label{blowup_eqn1}
\begin{cases}
y^a\partial_t U_i={\rm{div}}_X (y^a DU_i)\ \ {\rm{in}}\ \ P_2^+(0,0),\\
U_i(\partial \mathbb R^{m+1}_+\times (-4,4)) \subset {\epsilon_i}^{-1}\big(\mathbb S^{\ell-1}-U(0,1)\big), \\
\displaystyle y^a\frac{\partial U_i}{\partial y}\Big|_{y=0}\perp T_{U((r_i x, 0), r_i^2 t)} \mathbb S^{\ell-1}, 
\ (x,t)\in B_2^m(0)\times (-4,4).
\end{cases}
\end{equation} 
From \eqref{energy_ineq} we can deduce that
\begin{equation}\label{energy_ineq1}
\iint_{P_{1}^+(0,0)} |\partial_t U_i|^2\mu_adt\le C\iint_{P_{2}^+(0,0)} |DU_i|^2\mu_adt\le C.
\end{equation}
By the $y^a$-weighted Poincar\'e inequality, we have that
\begin{equation}\label{poincare}
\iint_{P_1^+(0,0)}  |U_i|^2\mu_adt\le C\iint_{P_{1}^+(0,0)} (|\partial_t U_i|^2+|DU_i|^2)\mu_adt\le C.
\end{equation} 
It follows from \eqref{decay6}, \eqref{energy_ineq1}, and \eqref{poincare} that we may assume that there exists
$U_\infty\in H^{1}(P_1^+(0,0), \mathbb R^{\ell})$
such that, after passing to a subsequence, 
$$(\partial_t U_i, DU_i) \rightharpoonup (\partial_t U_\infty, DU_\infty);
\ U_i\rightarrow U_\infty \ \  \ {\rm{in}}\ \ \ L^{2} (P_1^+(0,0), \mu_adt).
$$
Moreover, by using the $y^a$-weighted John-Nireneberg inequality \eqref{weighted-JN}, we can
further assume that
$$U_i\rightarrow U_\infty\ \ {\rm{in}}\ \ L^4(P_1^+(0,0), \mu_adt).$$  

It is not hard to check that 
\begin{equation}\label{decay9}
\iint_{P_{1}^+(0,0)}(|\partial_t U_\infty|^2+|DU_\infty|^2)\mu_adt\le 1, 
\end{equation}
and
\begin{equation}\label{blowup_limit_eqn1}
\begin{cases}
y^a\partial_t U_\infty={\rm{div}}_X(y^a DU_\infty)\ \ {\rm{in}}\ \ P_1^+(0,0),\\
U_\infty(B_1^m(0)\times (-1,1)) \subset T_{U(0,1)}\mathbb S^{\ell-1}, \\
\displaystyle y^a\frac{\partial U_\infty}{\partial y}\Big|_{y=0}\perp T_{U(0,1))} \mathbb S^{\ell-1}, 
\ (x,t)\in B_1^m(0)\times (-1,1).
\end{cases}
\end{equation} 
By choosing a suitable coordinate system in $\mathbb R^{\ell}$, we may assume $U(0,1)=(0',1)$ and $T_{U(0,1)}\mathbb S^{\ell-1}=\mathbb R^{\ell-1}\times \{0\}$ so that 
\begin{equation}\label{bc}
U_\infty^\ell=0, \ \  \displaystyle y^a\frac{\partial U_\infty^i}{\partial y}\Big|_{y=0}=0\ (1\le i\le \ell-1)\ \ {\rm{on}}\ \ B_1^m(0)\times (-1,1).
\end{equation}
Let us extend $U_\infty$ to $P_1(0,0)$ in such a way that $U_\infty^{\ell}$ is odd and $U_\infty^i, 1\le i\le \ell-1$, is even with respect to $y$. Then
it is easy to verify that $U_\infty$ solves 
\begin{equation}\label{blowuo_limt_eqn2}
|y|^a\partial_t U_\infty={\rm{div}}_X (|y|^a DU_\infty)\ \ {\rm{in}}\ \ P_1(0,0).
\end{equation}
Hence by the regularity theorem  that $U_\infty\in C^{\alpha_s}(P_1(0,0))$ and, for any $\theta\in (0,\frac14)$,
\begin{equation}\label{decay10}
\theta^{-(m+3+a)}\iint_{P_{\theta}^+(0,0)}|U_\infty-(U_\infty)^a_{(0,0),\theta, +}|^2\mu_adt\le C\theta^{2\alpha_s},
\end{equation}
where 
$$
\alpha_s=\begin{cases} {2s} & 0<s<\frac12,\\ 1 & \frac12\le s<1. \end{cases}
$$

With \eqref{decay10} at the hand,  we can apply \eqref{hybrid2} to reach a contradiction to \eqref{decay7} as follows.  
First of all, by the energy monotonicity inequality \eqref{local_monotonicity1} and the local energy inequality \eqref{local_energy_ineq1} for $U$,
we have that for any $z\in P_{r_i}^+(0,1)$ and $0<r<\theta_0{r_i}$, 
\begin{equation}\label{mono5'}
r^{-(m+1+a)}\iint_{P_r^+(z)}|DU|^2\mu_adt\le Cr_i^{-(m+1+a)}\iint_{P_{2r_i}^+(0,1)} |DU|^2\mu_adt,
\end{equation}
and
\begin{equation}\label{mono6}
r^{-(m+1+a)}\iint_{P_r^+(z)}r^2|\partial_tU|^2\mu_adt\le Cr_i^{-(m+1+a)}\iint_{P_{2r_i}^+(0,1)}|DU|^2\mu_adt.
\end{equation}
It follows from \eqref{mono5'} and \eqref{mono6} that 
\begin{align}\label{morrey7}
&\Big\|DU\Big\|_{M^{2,2}(P_{r_i}^+(0,1),\ \mu_adt)}+\Big\|\partial_tU\Big\|_{M^{2,4}(P_{r_i}^+(0,1), \  \mu_adt)}\nonumber\\
&\qquad\qquad\le C \Big(r_i^{-(m+1+a)}\iint_{P_{2r_i}^+(0,1)} |DU|^2\mu_adt\Big)^\frac12
\le C\varepsilon_i,
\end{align}
Let $\varepsilon_0>0$ be  the small constant given by Theorem \ref{hybrid1}.
We may assume $C\varepsilon_i\le \varepsilon_0$ for all $i$,  as $\varepsilon_i\to 0$.
Hence we can apply the inequality \eqref{hybrid2} of Theorem \ref{hybrid1}, with rescaling, to obtain that for any $0<\theta<\theta_0$, 
\begin{align}\label{hybrid3}
\theta^{-(m+1+a)}\iint_{P^+_\theta} |DU_i|^2\mu_adt&\le \varepsilon_0(2\theta)^{-(m+1+a)}\iint_{P^+_{2\theta}} ((2\theta)^2|\partial_t U_i|^2+|DU_i|^2)\mu_adt\nonumber\\
&+\varepsilon_0 (2\theta)^{-(m-1+a)}\sup_{t\in [-4\theta^2, 4\theta^2]}\int_{B_{2\theta}^+}|DU_i|^2\mu_a\nonumber\\ 
&+C\varepsilon_0^{-1}(2\theta)^{-(m+3+a)}\iint_{P^+_{2\theta}}|U_i-(U_i)^a_{0,2\theta,+}|^2\mu_adt\nonumber\\
&+C\varepsilon_0^{-1}(2\theta)^{-(m+3+a)}\iint_{P^+_{2\theta}}|U_i-(U_i)^a_{0,2\theta,+}|^4\mu_adt\nonumber\\
&+C\varepsilon_0^{-1} \varepsilon_i^{-2}(2\theta)^{-(m+3+a)}
\iint_{P^+_{2\theta}}(|\overline{U}_i|^2-1)^2\mu_adt,
\end{align}
where $\overline{U}_i(x,t)=U(r_ix, 1+r_i^2 t)$. 

Since $|\overline{U}_i|=1$ on $\widehat{P}_1(0,0)$, here $\widehat{P}_r(0,0)=B_r^m(0)\times (-r^2, r^2)$, it follows from the triangle inequality
that
\begin{align*}
&(2\theta)^{-(m+3+a)}\iint_{P^+_{2\theta}}(|\overline{U}_i|^2-1)^2\mu_adt
\le \Big(1+\|\overline{U}_i\|_{L^\infty(P_\frac12^+)}\Big)^2(2\theta)^{-(m+3+a)}\iint_{P^+_{2\theta}}|\overline{U}_i -\bar{c}_i(\theta)|^2\mu_adt\\
&\le C\Big(1+\|\overline{U}_i\|_{L^\infty(P_\frac12^+)}\Big)^2(2\theta)^{-(m+3+a)}\Big(\iint_{P^+_{2\theta}}|\overline{U}_i -(\overline{U}_i)^a_{0, 2\theta,+}|^2\mu_adt
+\iint_{P_{2\theta}^+} |\overline{U}_i-\bar{c}_i(\theta)|^2\mu_adt\Big)
\end{align*} 
where $\bar{c}_i(\theta)$ denotes the average of $\overline{U}_i$ over $\widehat{P}_{2\theta}(0,0)$ with respect to $dxdt$. 

From \eqref{morrey7}, and \eqref{dist-est} of Lemma \ref{dist-est0}, we have
\[
\big\|\overline{U}_i\big\|_{L^\infty(P_\frac12^+)}\le 1+C\varepsilon_i\le 2.
\]
Hence, after rescaling, we obtain 
\begin{align}\label{morrey8}
&\varepsilon_i^{-2}(2\theta)^{-(m+3+a)}\iint_{P^+_{2\theta}}(|\overline{U}_i|^2-1)^2\mu_adt\nonumber\\
&\le C(2\theta)^{-(m+3+a)}\Big(\iint_{P^+_{2\theta}}|{U}_i -({U}_i)^a_{0, 2\theta,+}|^2\mu_adt
+\iint_{P_{2\theta}^+}|{U}_i-{c}_i(\theta)|^2\mu_adt\Big),
\end{align} 
with ${c}_i(\theta)$ denotes the average of ${U}_i$ over $\widehat{P}_{2\theta}(0,0)$.

Since $U_i\rightarrow U_\infty$ in $L^4(P_{\theta_0}^+(0,0), \mu_adt)$ and $U_\infty$ satisfies \eqref{decay10}, we can estimate
\begin{align*}
&(2\theta)^{-(m+3+a)}\Big(\iint_{P^+_{2\theta}}|{U}_i -({U}_i)^a_{0, 2\theta,+}|^2\mu_adt
+\iint_{P_{2\theta}^+} |{U}_i-{c}_i(\theta)|^2\mu_adt\Big)\\
&\le 2(2\theta)^{-(m+3+a)}\Big(\iint_{P^+_{2\theta}}|{U}_\infty-({U}_\infty)^a_{0, 2\theta,+}|^2\mu_adt
+\iint_{P_{2\theta}^+} |{U}_\infty-{c}_\infty(\theta)|^2\mu_adt\Big)\\
&+4(2\theta)^{-(m+3+a)}\iint_{P^+_{2\theta}}|{U}_i-{U}_\infty|^2\mu_adt\\
&\le C\theta^2+C\theta^{-(m+3+a)} o_i(1). 
\end{align*} 
where $o_i(1)\to 0$ as $i\to \infty$.  Substituting this estimate into \eqref{morrey8}, we have
\begin{align}\label{morrey9}
\varepsilon_i^{-2}(2\theta)^{-(m+3+a)}\iint_{P^+_{2\theta}}(|\overline{U}_i|^2-1)^2\mu_adt
&\le C\theta^2+C\theta^{-(m+3+a)} o_i(1).
\end{align} 
Substituting \eqref{morrey9} and the scaled version of estimate \eqref{morrey7} into \eqref{hybrid3}, we arrive at 
\begin{align}\label{hybrid4}
\theta^{-(m+1+a)}\iint_{P^+_\theta} |DU_i|^2\mu_adt&\le C\varepsilon_0+C\theta^2+C\theta^{-m-3-a}o_i(1).
\end{align}
We now choose the parameters orderly: (i) First, choose $\varepsilon_0>0$ so small that $C\varepsilon_0<\frac18$,
(ii) Second, choose $\theta>0$  so small that $C\theta^2<\frac18$, and (iii) Finally, since  $o_i(1)\to 0$ as $i\to\infty$, 
choosing $i$ sufficiently large so that $C\theta^{-m-3-a}o_i(1)<\frac18$. This yields
\[
\theta^{-(m+1+a)}\iint_{P^+_\theta} |DU_i|^2<\frac38,
\]
which contradicts to \eqref{decay7}.
The proof is complete. \qed

\medskip
Iterations of \eqref{decay3} give rise to H\"older continuity of $U$. More precisely, we have
\begin{lemma}\label{holder-cont1} There exist $\epsilon_0>0$ and $\alpha_0\in (0,1)$ such that if $U$ is a stationary 
solution of heat flow of $s$-harmonic maps in $\R^{m+1}_+$,
which satisfies, for $(X_0,t_0)\in\partial\mathbb R^{m+1}_+\times (0,\infty)$ and $0<r<\frac{\sqrt{t_0}}{4}$,
\begin{equation}\label{small_cond}
(2r)^{-(m+1+a)}\iint_{P_{2r}^+(X_0,t_0)} |DU|^2 \mu_adt\leq \epsilon_0^2,
\end{equation}
then $U\in C^{\alpha_0,\frac{\alpha_0}2}_{X,t}(P_{\frac{r}2}^+(X_0,t_0))$, and,  for any $Z\in P_{r}^+(Z_0)$ and $0<\tau<r$,
\begin{equation}\label{holder-cont2}
\tau^{-(m+1+a)}\iint_{P_{\tau}^+(Z_1)}|DU|^2 \mu_adt\le \big(\frac{\tau}{r})^{2\alpha_0}  r^{-(m+1+a)} \iint_{P_{r}^+(Z_1)} |DU|^2\mu_a dt.
\end{equation}
\end{lemma}
\pf Let $Z_0=(X_0,t_0)$.  It follows from \eqref{small_cond} that 
\begin{equation}\label{small_cond1}
r^{-(m+1+a)}\iint_{P_r^+(X_*,t_*)} |DU|^2 \mu_adt\leq 2^{m+1+a}\epsilon_0^2
\end{equation}
holds for all $(X_*,t_*)\in \widehat{P}_r(Z_0)$. Hence we can repeatedly apply Lemma \ref{decay1} to conclude that there exists $\theta_0\in (0,1)$
such that for any $k\ge 1$,
\begin{equation}\label{decay-k}
(\theta_0^kr)^{-(m+1+a)}\iint_{P_{\theta_0^kr}^+(X_*,t_*)} |DU|^2\mu_adt\le 2^{-k} r^{-(m+1+a)} \iint_{P_{r}^+(X_*,t_*)} |DU|^2\mu_adt
\end{equation}
holds for all $(X_*,t_*)\in \widehat{P}_r(Z_0)$. This, combined with the local energy inequality \eqref{energy_ineq}, yields
\begin{align}\label{decay-kk}
&(\theta_0^kr)^{-(m+1+a)}\iint_{P_{\theta_0^kr}^+(X_*,t_*)}\big(|DU|^2+(\theta_0^k r)^2|\partial_t U|^2\big)\mu_adt\nonumber\\
&\qquad\qquad\qquad\le C2^{-k} r^{-(m+1+a)} \iint_{P_{r}^+(X_*,t_*)} |DU|^2\mu_adt
\end{align}
holds for all $(X_*,t_*)\in \widehat{P}_r(Z_0)$. Standard interpolation then implies that there exists $\alpha_1\in (0,1)$ such that
\begin{align}\label{decay-kkk}
&\tau^{-(m+1+a)}\iint_{P_{\tau}^+(X_*,t_*)}\big(|DU|^2+\tau^2|\partial_t U|^2\big)\mu_adt\nonumber\\
&\qquad\qquad\qquad\le C\big(\frac{\tau}{r}\big)^{2\alpha_1} r^{-(m+1+a)} \iint_{P_{r}^+(X_*,t_*)} |DU|^2\mu_adt
\end{align}
holds for all $Z_*=(X_*,t_*)\in \widehat{P}_r(Z_0)$ and $0<\tau<r$.

On the other hand, by the interior H\"older continuity of $U$, we have that there exists $\alpha_2\in (0,1)$ such that 
for any $Z_*=(x_*,y_*,t_*)\in P_{r}^+(Z_0)$ with $y_*>0$, 
\begin{align}\label{decay-kkkk}
&\tau^{-(m+1+a)}\iint_{P_{\tau}(Z_*)}\big(|DU|^2+\tau^2|\partial_t U|^2\big)\mu_adt\nonumber\\
&\qquad\qquad\le C(\frac{\tau}r)^{2\alpha_2} y_*^{-(m+1+a)} \iint_{P_{y_*}^+(Z_*)} |DU|^2\mu_adt\nonumber\\
&\qquad\qquad\le C(\frac{\tau}r)^{2\alpha_2} (2y_*)^{-(m+1+a)} \iint_{P_{2y_*}^+((X_*,0,t_*))} |DU|^2\mu_adt
\end{align}
holds for all $0<\tau<y_*$.

Combining \eqref{decay-kkk} and \eqref{decay-kkkk}, it is not hard to show that
\begin{equation}
\tau^{-(m+1+a)} \int_{P_\tau^+(Z_*)} \big(|DU|^2+\tau^2|\partial_t U|^2)\mu_adt\le C(\frac{\tau}r)^{2\alpha_0} r \iint_{P_{r}^+(Z_0)} |DU|^2\mu_adt
\end{equation}
holds for all $Z_*\in P_r^+(Z_0)$ and $0<\tau<r$, where $\alpha_0=\min\{\alpha_1,\alpha_2\}$. 
Thus, by Morrey's decay Lemma (see \cite{Giaquinta}) and Lemma \ref{BMO_a} and Lemma \ref{weighted-poincare} below, we conclude that
$U\in C^{\alpha_0, \frac{\alpha_0}2}_{X,t}(P_r^+(Z_0))$ and \eqref{holder-cont2} holds.
\qed

\section{Proof of Theorem \ref{partial_reg}: Higher order regularity}

This section is devoted to the proof of $C^{1,\alpha}$-regularity of the heat flow of $s$-harmonic maps $u=U\big|_{\partial\R^{m+1}_+\times (0,\infty)}$
in a neighborhood of  $Z_0=(X_0,t_0)\in \partial\R^{m+1}_+\times (0,\infty)$ where $U$ is H\"older continuous. 
This amounts to show that $(\partial_t ^ku, D_x^ku)$ is  continuous near $Z_0$ for all $k\ge 1$. 

To achieve this, we first choose a small $R_0>0$ and extend $U$ from $P_{R_0}^+(Z_0)$ to $P_{R_0}(Z_0)$ such that the extension $\overline{U}$
solves a $|y|^a$-weighted parabolic system of approximated harmonic maps. 

Since $U$ is H\"older continuous near $Z_0=(X_0,t_0)\in \partial\R^{m+1}_+\times (0,\infty)$, there exists $R_0>0$ such that
\begin{equation}\label{bound_U}
\frac12\le |U(Z)|\le 1, \ \forall Z\in P_{R_0}^+(Z_0).
\end{equation}
Now we define the extension map $\overline{U}:P_{R_0}(Z_0)\to \R^{\ell}$ by 
\begin{equation*} \overline{U}(X,t)=
\begin{cases} U(X,t), & (X,t)=((x,y),t)\in P_{R_0}^+(Z_0),\\
\displaystyle\big(2\frac{U}{|U|}-U\big)(x,-y,t), & (X,t)=((x,y),t)\in P_{R_0}^-(Z_0).
\end{cases}
\end{equation*}
We remark that this type of extension has played crucial roles in the study of boundary regularity for heat flow of harmonic maps 
subject to free boundary conditions (see for example \cite{Struwe1, ChenLin1}). 

\medskip
By direct calculations and application of the boundary condition of $U$, we obtain
\begin{lemma} The extension map $\overline{U}: P_{R_0}(Z_0)\to\R^{\ell}$ solves
\begin{equation}\label{ext_eqn}
|y|^a \partial_t\overline{U}-{\rm{div}}_X\big(|y|^a D\overline{U})=|y|^a\Gamma,
\end{equation}
where
\begin{equation}\label{ext_rhs}
\Gamma(X,t)=\begin{cases} 0, & (X,t)\in P_{R_0}^+(Z_0),\\
2\Big(\displaystyle\frac{|DU|^2U}{|U|^3}+2\frac{DU\langle DU, U\rangle}{|U|^3}-3\frac{|\langle DU, U\rangle|^2 U}{|U|^5}\Big)(x,-y,t), & (X,t)\in P_{R_0}^-(Z_0).
\end{cases}
\end{equation}
\end{lemma}

Observe that \eqref{ext_eqn} and \eqref{ext_rhs} can be viewed
as a $|y|^a$ weighted parabolic harmonic map like system, since
\[
|\Gamma(X,t)|\le \begin{cases} 0, & (X,t)\in P_{R_0}^+(Z_0),\\
C |DU|^2(x,-y,t), & (X,t)\in P_{R_0}^-(Z_0).
\end{cases}
\]

\bigskip
It turns out that some of the techniques to prove higher order boundary regularity for the heat flow of harmonic maps
can be adapted to establish the smoothness of $u$ near $Z_0$.  We divide the proof of smoothness into three steps.

\medskip
\noindent{\it Step 1}. Improving the H\"older continuity of $U$. We want to show $\overline{U}\in C^{\gamma_a, \frac{\gamma_a}2}_{X,t}(P_{\frac{R_0}2}(Z_0))$, where
\[
\gamma_a=\begin{cases} 2s, &\  0<s<\frac12,\\
1, & \ \frac12\le s<1. \end{cases}
\]
To show this, let $R\in (\frac{R_0}2, R_0)$ and $V: P_{R}(Z_0)\to \R^{\ell}$ solve 
\begin{equation}\label{auxi-eqn}
\begin{cases}
|y|^a\partial_t V-{\rm{div}}_X(|y|^a DV)=0&\ \ {\rm{in}} \  \ P_{R}(Z_0),\\
V=\overline{U}&\ \ {\rm{on}}\ \ \partial_p P_R(Z_0).
\end{cases}
\end{equation}
Here $\partial_p P_R(Z_0)$ denotes the parabolic boundary of $P_R(Z_0)$. 

Since $\overline{U}$ is H\"older continuous
with exponent $\alpha_0\in (0,1)$,  the maximum principle implies that 
\[
{\rm{osc}}_{P_R(Z_0)} V\le {\rm{osc}}_{P_R(Z_0)} \overline{U}\le CR^{\alpha_0}.
\] 
Testing \eqref{ext_eqn} and \eqref{auxi-eqn} with $\overline{U}-V$ and integrating the resulting equations over $P_{R}(Z_0)$, 
then subtracting each other, we obtain
\begin{align}\label{decay1'}
\iint_{P_R(Z_0)} |D(\overline{U}-V)|^2|y|^adXdt&\le C\|\overline{U}-V\|_{L^\infty(P_R(Z_0))} \iint_{P_R^+(Z_0)} |DU|^2|y|^adXdt\nonumber\\
&\le C\big({\rm{osc}}_{P_R(Z_0)} \overline{U}\big) R^{m+1+a+2\alpha_0}\nonumber\\
&\le CR^{m+1+a+3\alpha_0}.
\end{align} 
On the other hand, it is well-known (see \cite{DongPhan}) that $V\in C^{\gamma_a,\frac{\gamma_a}2}(P_{\frac{3R}4}(Z_0))$, and
\begin{align}\label{decay2'}
(\theta R)^{-(m+1+a)}\iint_{P_{\theta R}(Z_0)} |DV|^2|y|^adXdt&\le C \theta^{2\gamma_a} R^{-(m+1+a)}\iint_{P_R(Z_0)} |DV|^2|y|^adXdt\nonumber\\
&\le C\theta^{2\gamma_a} R^{-(m+1+a)}\iint_{P_R(Z_0)} |D\overline{U}|^2|y|^adXdt\nonumber\\
&\le C\theta^{2\gamma_a} R^{-(m+1+a)}\iint_{P_R^+(Z_0)} |D{U}|^2\mu_adt.
\end{align}
holds for all $0<\theta<\frac34$.

Combining \eqref{decay1'} with \eqref{decay2'}, we obtain that for any $\theta\in (0,\frac34)$, there holds
\begin{align*}
&(\theta R)^{-(m+1+a)}\iint_{P_{\theta R}^+(Z_0)} |DU|^2\mu_adt\\
&\le 2(\theta R)^{-(m+1+a)}\iint_{P_{\theta R}(Z_0)} \Big(|DV|^2+|D(\overline{U}-V)|^2\Big)|y|^adXdt\\
&\le C\theta^{2\gamma_a} R^{-(m+1+a)} \iint_{P_R^+(Z_0)} |DU|^2\mu_adt+C\theta^{-(m+1+a)} R^{3\alpha_0}.
\end{align*}

Without loss of generality, we start the iterations by assuming $3\alpha_0\le 2\gamma_a$.  Then, by Giaquinta's iteration Lemma
\cite{Giaquinta}, we conclude that for $0<\rho<R\le \frac{R_0}2$, 
\begin{equation}\label{hold1}
\rho^{-(m+1+a)}\iint_{P_{\rho}^+(Z_0)} |DU|^2\mu_adt\le C\Big[(\frac{\rho}{R})^{3\alpha_0} R^{-(m+1+a)}\iint_{P_{R}^+(Z_0)} |DU|^2\mu_adt+\rho^{3\alpha_0}\Big].
\end{equation}
Together with the local energy inequality \eqref{local_energy_ineq1}, this leads to 
\begin{align}\label{hold1'}
&\rho^{-(m+1+a)}\iint_{P_{\rho}^+(Z_0)} (\rho^2|\partial_t U|^2+|DU|^2)\mu_adt\nonumber\\
&\qquad\qquad\qquad\le C\Big[(\frac{\rho}{R})^{3\alpha_0} R^{-(m+1+a)}\iint_{P_{R}^+(Z_0)} |DU|^2\mu_adt+\rho^{3\alpha_0}\Big].
\end{align}

On the other hand, if $Z_*=(x_*, y_*, t_*)\in \R^{m+1}_+\times (0,\infty)$, then by the interior smoothness of $U$, it is not hard to show that for $0<\rho<y_*$, we have
\begin{align}\label{hold2}
&\rho^{-(m+1+a)}\iint_{P_{\rho}(Z_*)} (\rho^2|\partial_t U|^2+|DU|^2)\mu_adt\nonumber\\
&\le C\big(\frac{\rho}{y_*}\big)^{2\gamma_a} y_*^{-(m+1+a)}\iint_{P_{y_0}(Z_*)} |DU|^2\mu_adt\nonumber\\
&\le C\big(\frac{\rho}{y_*}\big)^{2\gamma_a} y_*^{-(m+1+a)}\iint_{P_{2y_*}^+(x_*, 0, t_*)} |DU|^2\mu_adt.
\end{align} 
Based on \eqref{hold1} and \eqref{hold2}, we can show that the following inequality
\begin{align}\label{hold3}
&\rho^{-(m+1+a)}\iint_{P^+_{\rho}(Z_1)} (\rho^2|\partial_t U|^2+|DU|^2)\mu_adt\nonumber\\
&\qquad\qquad\qquad\qquad\qquad\le C\Big[\big(\frac{\rho}{R_0}\big)^{2\alpha_1} R_0^{-(m+1+a)}\iint_{P_{R_0}^+(Z_0)} |DU|^2\mu_adt
+\rho^{2\alpha_1}\Big]
\end{align} 
holds for all $Z_1\in P^+_{\frac{R_0}2}(Z_0)$ and $0<\rho\le \frac{R_0}4$, where $\alpha_1=\min\{\gamma_a, \frac{3\alpha_0}2\}>\alpha_0$.

Applying  the $y^a$-weighted Morrey decay Lemma \cite{Giaquinta} to \eqref{hold3}, we have $U\in C^{\alpha_1, \frac{\alpha_1}2}_{X,t}(\overline{P_{\frac{R_0}2}^+(Z_0)})$.
If $\alpha_1<\gamma_a$, then we can repeat the same argument to improve $U\in C^{\alpha_2,\frac{\alpha_2}2}_{X,t}(\overline{P_{\frac{R_0}4}^+(Z_0)})$,
with $\alpha_2=\min\{\gamma_a, \frac{9\alpha_0}4\}.$ It is clear that it takes at most $k$-times for some $k\ge 1$ so that the H\"older exponent 
$\alpha_k=\gamma_a$ and $U\in C^{\gamma_a, \frac{\gamma_a}2}_{X,t}(\overline{P^+_{\frac{R_0}{2^k}}(Z_0)})$.

\medskip
\noindent{\it Step 2}.  Regularity of $\lambda_u$ for $0<s<\frac12$. Since $\frac12\le s<1$ can be handled similarly,
we focus on the case $0 <s<\frac12$. From Step 1 and suitable scaling, we may assume that $Z_0=(0, 1)\in \partial \mathbb{R}^{m+1}_+\times (0,\infty)$
and $u\in C^{2s, s}_{X,t}\big(\widehat{P}_2^+(0,1)\big)$. We want to estabish the following:

\medskip
\noindent{\bf Claim}. {\it For any small $0<\varepsilon<s$,  $\lambda_u\in C^{2(s-\varepsilon), s-\varepsilon}_{x,t}\big(\widehat{P}_\frac14^+(0,1)\big)$}.
To see this, let's recall
the parabolic distance on $\R^m\times (0,\infty)$ by
\[
d((x,t), (y,\tau))=|x-y|+\sqrt{|t-\tau|}, \ \ (x,t),\ (y,\tau)\in \R^m\times (0,\infty).
\]
Given two different points $z_1=(x_1, t_1), z_2=(x_2, t_2)\in \widehat{P}_\frac14^+(0,1) $. Then $0<\delta=d\big(z_1, z_2\big)\le \frac12$, and
\begin{align*}
&\big|\lambda_u(z_1)-\lambda_u(z_2)\big|\\
&\le\Big|\iint_{\{d((y,\tau), z_1)<2\delta\}} \big[|u(z_1)-u(y,\tau)|^2K_s(x_1-y,t_1-\tau)-|u(z_2)-u(y,\tau)|^2K_s(x_2-y, t_2-\tau)\big]\Big|\\
&+\Big|\iint_{\{d((y,\tau), z_1)\ge 2\delta\}} \big[|u(z_1)-u(y,\tau)|^2K_s(x_1-y,t_1-\tau)-|u(z_2)-u(y,\tau)|^2K_s(x_2-y, t_2-\tau)\big]\Big|\\
&=I+II.
\end{align*}
Since 
$$|u(z)-u(w)|\le Cd(z,w)^{2s}, \ \forall z, w\in \widehat{P}_2^+(0,1),$$
and
$$
|K_s(x-y,t-\tau)|\le \frac{C}{d((x,t), (y,\tau))^{n+2+2s}}, \ \ (x,t), \ (y,\tau)\in\R^m\times (0,\infty),
$$
and $d(z,z_1)\le 2\delta$ implies $d(z, z_2)\le 3\delta$ by the triangle inequality, we can bound $I$ by
\begin{align*}
|I|&\le \iint_{\{d((y,\tau), z_1)<2\delta\}} \big|u(z_1)-u(y,\tau)|^2|K_s(x_1-y,t_1-\tau)|\\
&\quad+\iint_{\{d((y,\tau), z_2)<3\delta\}} |u(z_2)-u(y,\tau)|^2|K_s(x_2-y, t_2-\tau)|\\
&\le C\iint_{\{d((y,\tau), z_1)<2\delta\}}\frac{1}{d((y,\tau), z_1)^{n+2-2s}}+\iint_{\{d((y,\tau), z_2)<3\delta\}}\frac{1}{d((y,\tau), z_2)^{n+2-2s}}\\
&\le C\int_0^{2\delta} r^{-1+2s}\,dr+C\int_0^{3\delta} r^{-1+2s}\,dr\le C\delta^{2s}.
\end{align*} 
For $II$, we have
\begin{align*}
|II|&\le \iint_{\{d((y,\tau), z_1)\ge 2\delta\}} |u(z_1)-u(y,\tau)|^2\Big|K_s(x_1-y,t_1-\tau)-K_s(x_2-y, t_2-\tau)\Big|\\
&+\iint_{\{d((y,\tau), z_1)\ge 2\delta\}} \Big||u(z_1)-u(y,\tau)|^2-|u(z_2)-u(y,\tau)|^2\Big|K_s(x_2-y, t_2-\tau)\big|\\
&=III+IV.
\end{align*}
Since 
\[
\Big|K_s(x_1-y,t_1-\tau)-K_s(x_2-y, t_2-\tau)\Big|\le \frac{Cd(z_1,z_2)}{d((x_1,t_1), (y,\tau))^{n+3+2s}}, \  \  {\rm{for}}\ \ d((y,\tau),z_1)\ge 2\delta,
\]
and $|u(z)-u(w)|^2\le 4$ for all $z,w\in \R^m\times (0,\infty)$, we can bound $III$ by
\begin{align*}
\big|III\big|&\le  \iint_{\{2\delta\le d((y,\tau), z_1)<2\}}|u(z_1)-u(y,\tau)|^2\Big|K_s(x_1-y,t_1-\tau)-K_s(x_2-y, t_2-\tau)\Big|\\
&+\iint_{\{d((y,\tau), z_1)\ge 2\}}|u(z_1)-u(y,\tau)|^2\Big|K_s(x_1-y,t_1-\tau)-K_s(x_2-y, t_2-\tau)\Big|\\
&\le C\delta\Big(\iint_{\{2\delta\le d((y,\tau), z_1)<2\}}\frac{d((y,\tau),z_1)^{4s}}{d((y,\tau),z_1)^{n+3+2s}}+\iint_{\{d((y,\tau), z_1)\ge 2\}}\frac{4}{d((y,\tau),z_1)^{n+3+2s}}\Big)\\
&\le C\delta\int_{2\delta}^2 r^{-(n+3-2s)} r^{n+1}\,dr+ C\delta\int_2^\infty r^{-(n+3+2s)} r^{n+1}\,dr\\
&\le C\delta^{2s}+C\delta\le C\delta^{2s},
\end{align*}
as $0<2s<1$. 

For $IV$, since $d(w,z_1)\ge 2\delta$ implies $d(w,z_2)\ge\delta$ by the triangle inequality, we have
\begin{align*}
\big|IV\big|&\le C\iint_{\{d((y,\tau), z_2)\ge \delta\}} \big|\langle u(z_1)-u(z_2), u(y,\tau)\rangle\big| d((y,\tau), z_2)^{-(n+2+2s)}\\
&\le C\iint_{\{d((y,\tau), z_2)\ge \delta\}} \frac{\big|\langle u(z_1)-u(z_2), u(y,\tau)-u(z_2)\rangle\big| +|u(z_1)-u(z_2)|^2}{d((y,\tau), z_2)^{n+2+2s}}\\
&\le C\iint_{\{d((y,\tau), z_2)\ge \delta\}} \frac{|u(z_1)-u(z_2)||u(y,\tau)-u(z_2)| +|u(z_1)-u(z_2)|^2}{d((y,\tau), z_2)^{n+2+2s}}\\
&\le C\delta^{2s}\iint_{\{\delta\le d((y,\tau), z_2)\le 1\}} \frac{d((y,\tau),z_2)^{2s}}{d((y,\tau), z_2)^{n+2+2s}}
+C\delta^{2s}\iint_{\{d((y,\tau), z_2)\ge 1\}} \frac{4}{d((y,\tau), z_2)^{n+2+2s}}\\
&+C\delta^{4s} \iint_{\{d((y,\tau), z_2)\ge \delta\}}\frac{1}{d((y,\tau), z_2)^{n+2+2s}}\\
&\le C\delta^{2s}\Big(\int_\delta^1 r^{-(n+2+2s)}r^{2s} r^{n+1}\,dr+\int_1^\infty r^{-(n+2+2s)} r^{n+1}\,dr
+\delta^{2s} \int_\delta^\infty r^{-(n+2+2s)} r^{n+1}\,dr\Big)\\
&\le C\delta^{2s}(1+|\ln\delta|). 
\end{align*}
Putting all these estimates together, we obtain that
\[
|\lambda_u(z_1)-\lambda_u(z_2)|\le Cd(z_1, z_2)^{2s}\big(1+|\ln d(z_1, z_2)|\big), \ \ \forall z_1, z_2\in \widehat{P}_\frac14^+(Z_0).
\]
This clearly yields the claim.

\medskip
\noindent{\it Step 3}. Smoothness of $u=U\big|_{\partial\R^{m+1}_+\times (0,\infty)}$ via Schauder estimate.  
For $Z_0=(x_0, 0, t_0)\in \partial\R^{m+1}_+\times (0,\infty)$,
let $0<R_0<\frac{\sqrt{t_0}}2$ be such that $u\in C^{\gamma_a, \frac{\gamma_a}2}_{x,t}\big(D_{R_0}(x_0)\times (t_0-R_0^2, t_0+R_0^2)\big)$. 
We want to show $u\in C^\infty\big(\widehat{P}_{\frac{R_0}{16}}(x_0,t_0)\big)$.  For simplicity, we only consider $0<s<\frac12$
and sketch the proof by assuming $(x_0,t_0)=(0, 1)$
and $R_0=2$. From Step 2, we know that $\lambda_u(x,t) u(x,t)\in C^{\gamma_a, \frac{\gamma_a}2}_{x,t}(\widehat{P}^+_\frac14(0,1))$.
Let $\eta\in C_0^\infty(\R^m\times (0,\infty))$ be a cut-off function of $\widehat{P}^+_\frac15(0,1)$ and let $v:\mathbb R^m\times \R\to\R^\ell$ solve
\[
(\partial_t-\Delta)^s v=\eta^2\lambda_u u\ \ {\rm{in}}\ \ \R^m\times \R.
\] 
Then, by the Schauder estimate for fractional parabolic equations (see \cite[Theorem 1.19]{ST}), we deduce that
$v\in C^{4s-2\varepsilon, 2s-\varepsilon}_{x,t}(\R^m\times (0,\infty))$, $\forall 0<\varepsilon<s$.

Since $u$ solves \eqref{eq:fhf1}, that is,  
\[
(\partial_t-\Delta)^s u=\lambda_u u, \ \ {\rm{in}}\ \  \R^m\times (0,\infty),
\]
If we denote $w=u-v$. Then 
\[
(\partial_t-\Delta)^s w=0\ \ {\rm{in}}\ \ \widehat{P}^+_\frac15(0,1).
\] 
Applying the Schauder estimate for fractional parabolic equations once again, we deduce that 
$w\in C^{\infty}(\widehat{P}^+_\frac15(0,1))$. Thus $u\in C^{4s-2\varepsilon, 2s-\varepsilon}_{x,t}(\widehat{P}^+_\frac15(0,1)))$, $\forall 0<\varepsilon<s$.

By the bootstrap argument, we can conclude that $u\in C^{2ks-2\varepsilon, ks-\varepsilon}_{x,t}(\widehat{P}^+_\frac18(0,1))$ for all $k\ge 1$
and $0<\varepsilon<s$. In particular, $u$ is smooth in $\widehat{P}^+\frac18(0,1)$. This establishes the smoothness of $u$ near $Z_0$.

\section{Appendix} 

In this section, we will provide several important Lemmas that are used in the previous sections.  We begin this section by pointing out  that
for $-1<a<1$, the weight function $$w_a(x,y)=|y|^a,  \ \ (x,y)\in \R^{m+1},$$
 is a Muckenhoupt $A_2$-weight function
 on $\R^{m+1}$,
namely,
\begin{equation}\label{A2}
\big[w_a\big]_{A_2}=\sup_{Q}\Big(\frac{1}{|Q|}\iint_Q w_a\,dXdt\Big)\Big(\frac{1}{|Q|}\iint_Q w_{-a}\,dXdt\Big)<+\infty,
\end{equation}
where the supremum is taken among all cubes $Q\subset\R^{m+1}$.

Now, we recall the definition of $y^a$-weighted BMO spaces in $\R^{m+1}_+\times (0,\infty)$. For any open set $O\subset \R^{m+1}\times (0,\infty)$,  
denote  $O_+=O\cap \overline{\R^{m+1}_+}\times (0,\infty)$ and define ${\rm{BMO}}(O_+,\mu_adt)$ to the space consisting of 
all $U\in L^1_{\rm{loc}}(O_+, \mu_adt)$ such that 
\begin{align}\label{weighted-bmo}
\big[U\big]^2_{{\rm{BMO}}(O_+,\ \mu_adt)}=\sup\Big\{\frac{1}{\displaystyle|P_r^+(z)|_a} \iint_{P_r^+(z)}|U-(U)^a_{P_r^+(z)}|^2\mu_adt:
P_r(z)\subset O_+\Big\}<\infty,
\end{align}
where 
$$|P_r^+(z)|_a=\iint_{P_r^+(z)} \mu_adt
\  {\rm{and}}\  (U)^a_{P_r^+(z)}=\frac{1}{|P_r^+(z)|_a} \iint_{P_r^+(z)} U \mu_adt.$$
Notice that when $a=0$, 
we have $\mu_0dt=dXdt$ and  ${\rm{BMO}}(O_+,\mu_adt)$ reduces to the standard BMO space ${\rm{BMO}}(O_+)$ with respect to the Lebesgue measure. 

The following simple fact on the relationship between the BMO space and the $y^a$-weighted BMO space plays an important role in this paper.
\begin{lemma} \label{BMO_a} If $U\in {\rm{BMO}}(O_+,\mu_adt)$, then $U\in BMO(O_+)$, and
\begin{equation}\label{bmo-bmoa} 
\big[U\big]_{{\rm{BMO}}(O_+)}\le C\big[U\big]_{{\rm{BMO}}(O_+,\  \mu_adt)}.
\end{equation}
\end{lemma}
\pf  For any $f\in L^1_{\rm{loc}}(O_+, \R^{\ell})$ and $P_r(z)\subset O_+$, by \eqref{A2} we have
\begin{align}\label{holdera}
\frac{1}{|P_r(z)|}\iint_{P_r(z)} |f| dXdt&\le \frac{1}{|P_r(z)|} \big(\iint_{P_r(z)} |f|^2 \mu_a dt\big)^\frac12 \big(\iint_{P_r(z)} y^{-a} dXdt\big)^\frac12\nonumber\\
&\le  C\big(\iint_{P_r(z)} |f|^2 \mu_a dt\big)^\frac12\big(\iint_{P_r(z)} y^a\,dXdt\big)^{-\frac12}\nonumber\\
&\le C\big(\frac{1}{|P_r(z)|_a}\iint_{P_r(z)} |f|^2 \mu_a dt\big)^\frac12.
\end{align} 
Applying \eqref{holdera} with $f=U-(U)^a_{P_r(z)}$, we can estimate
\begin{align*}
\inf_{c\in\R^{\ell}} \frac{1}{|P_r(z)|}\iint_{P_r(z)} |U-c| dXdt
&\le \frac{1}{|P_r(z)|}\iint_{P_r(z)} |U-(U)^a_{P_r(z)}|dXdt\\
&\le C\Big(\frac{1}{\displaystyle|P_r(z)|_a}\iint_{P_\tau(z)} |U-(U)^a_{P_r(z)}|^2\mu_adt\Big)^\frac12\\
&\le C \Big[U\Big]_{{\rm{BMO}}(O_+,\ \mu_adt)},
\end{align*}
this,  after taking supremum over all $P_r(z)\subset O_+$, implies \eqref{bmo-bmoa}.
\qed

\medskip
It is well-known that the following  $y^a$-weighted Poincar\'e inequality and Sobolev inequality hold, thanks to $y^a$ being a $A_2$-weight.

\begin{lemma}\label{weighted-poincare} There exists $C>0$ such that for any $f\in H^1_{\rm{loc}}(\R^{m+1}_+\times (0,\infty); \mu_adt)$, there holds
\begin{equation}\label{weight-poincare}
\iint_{P_R^+(Z)} |f-(f)^a_{P_R^+(z)}|^2\mu_adt\le CR^2\iint_{P_R^+(Z)} (|Df|^2+R^2|\partial_t f|^2)\mu_adt,
\end{equation}
and 
\begin{align}\label{weight-sobolev}
&\Big(\frac{1}{|P_R^+(Z)|_a}\iint_{P_R^+(Z)} |f-(f)^a_{P_R^+(z)}|^{2_a^*}\mu_adt\Big)^{\frac{1}{2^*}}\nonumber\\
&\qquad\qquad\qquad\qquad\le CR\Big(\frac{1}{|P_R^+(Z)|_a}\iint_{P_R^+(Z)} (|Df|^2+R^2|\partial_t f|^2)\mu_adt\Big)^\frac12,
\end{align}
for any $P_R(Z)\subset \R^{m+1}\times (0,\infty)$, where $2_a^*$ is the $y^a$-weighted Sobolev exponent:
\begin{equation}\label{sobolev-exponent}
2_a^*=\begin{cases} \frac{2(m+1-a)}{m-1-a}, &\ m\ge 2, \ {\rm{or}}\  m=1 \ {\rm{and}} \ a<0,\\
\infty, & \ m=1 \ {\rm{and}} \ a>0.
\end{cases}
\end{equation}
\end{lemma}

In the proof of Lemma \ref{dist-est0}, we need the following modified Poincar\'e inequality. 
Let $D_R(x_0)\subset\R^m$ be the ball with center $x_0\in\R^m$ and radius $R>0$ so that $D_R(x_0)=\{x\in \R^m:  |x-x_0|\le R\}$.
Denote $\widehat{P}_R(x_0,t_0)=D_R(x_0)\times [t_0-R^2, t_0+R^2]$.

\begin{lemma} \label{poincare1} There exists a constant $C>0$ such that for any $x_0\in\mathbb R^m$, $t_0>0$ and $0<R<\sqrt{t_0}$, it holds that
\begin{align}\label{poincare2}
&R^{-1}\iint_{P_R^+(x_0,t_0)}|U-[u]_{(x_0,t_0), R}|\,dXdt+\iint_{\widehat{P}_R(x_0,t_0)}|u-[u]_{(x_0,t_0), R}|\,dxdt\nonumber\\
&\le C \iint_{P_R^+(x_0,t_0)}
(|DU|+R|\partial_tU|)\,dXdt,
\end{align}
where $u=U\big|_{\partial R^{m+1}_+\times (0,\infty)}$, and $ [u]_{(x_0,t_0), R}$ denotes the average of $u$ over $\widehat{P}_R(x_0,t_0)$
with respect to $dxdt$.
\end{lemma}
\pf  Here we give a sketch of proof.  By scaling, it is sufficient to show \eqref{poincare2} for $R=1$. 
Suppose it were false. There would exist a sequence
$U^i$, with $(\partial_tU^i, DU^i)\in L^2(P_1^+(x_0,t_0))$, such that
\begin{align}\label{poincare3}
\alpha_i&:=\|U^i-[u^i]_{(x_0,t_0),1}\|_{L^1(P_1^+(x_0,t_0))}+\big\|u^i-[u^i]_{(x_0, t_0), 1}\big\|_{L^1(\widehat{P}_1(x_0,t_0))}\nonumber\\
&>i \iint_{P_1^+(x_0,t_0)}
(|DU^i|+|\partial_tU^i|)\,dXdt.
\end{align}
Replacing $U^i$ by
$$\displaystyle V^i=\frac{U^i-[u^i]_{(x_0, t_0),1}}{\alpha_i},$$
we may further assume that
$$[u^i]_{(x_0, t_0),1}=0 \ {\rm{and}}\  \big\|U^i-[u^i]_{(x_0,t_0), 1}\|_{L^1(P^+_1(x_0,t_0)}+\big\|u^i\big\|_{L^1(\widehat{P}_1(x_0, t_0)}=1.$$
Hence we have 
\begin{equation}\label{poincare4}
\iint_{P_1^+(x_0,t_0)}(|DU^i|+|\partial_tU^i|)\,dXdt\leq \frac{1}{i}.
\end{equation}
Recall by the standard Poincar\'e inequality, we have that
\begin{align}
\iint_{P_1^+(x_0,t_0)} |U^i|\,dXdt &\le C\big(\iint_{{P}^+_1(x_0,t_0)}(|DU^i|+|\partial_tU^i|)\,dXdt
+\iint_{\widehat{P}_1(x_0,t_0)} |u^i|\,dxdt\big)\nonumber\\
&\le C(1+i^{-1})\le 2C.
\end{align}
Thus $\{U^i\}\subset W^{1,1}(P_1^+(x_0,t_0))$ is a bounded sequence. Since 
$$W^{1,1}(P_1^+(x_0,t_0))\subset 
L^1(\widehat{P}_1(x_0,t_0))$$
is compact,  we can conclude that there exists a $U\in BV(P_1^+(x_0,t_0))$ such that
$$U^i\rightharpoonup U\ \ {\rm{ in }}\ \ BV(P_1^+(x_0,t_0)), \ \ \ U_i\rightarrow U\ \ {\rm{  in }}\  \ L^1(P_1^+(x_0,t_0)),$$
and 
$$u^i\rightarrow u=U\big|_{\widehat{P}_1(x_0,t_0)} \ \ {\rm{ in }}\ \ L^1(\widehat{P}_1(x_0,t_0)).$$
Hence 
$$[u]_{(x_0, t_0), 1}=0,  \ \|U\|_{L^1(P_1^+(x_0,t_0))}+\|u\|_{L^1(\widehat{P}_1(x_0,t_0))}=1.$$
On the other hand, it follows from \eqref{poincare4} and the lower semicontinuity that  
$$(\partial_tU, DU)=(0,0)\ \ \ {\rm{in}}\ \ \ P_1^+(x_0,t_0).$$
Hence $U$ is a constant, which, combined with $[u]_{(x_0, t_0), 1}=0$, implies that $U\equiv 0$. This contradicts to
$\|U\|_{L^1(P_1^+(x_0,t_0))}+\|u\|_{L^1(\widehat{P}_1(x_0,t_0))}=1$.  This completes the proof. \qed

\bigskip
For the proof of Theorem \ref{hybrid1}, we need the John-Nirenberg inequality in ${\rm{BMO}}(O_+,\mu_adt)$, whose proof can be found
in Cruz-Uribe--Martell--P\'erez \cite{CMP}. 

\begin{lemma}\label{weight-JN} Let $f\in {\rm{BMO}}(P_\frac12^+(0,1),\ \mu_adt)$. Then for any $1\le p<\infty$, there exists $C_p>0$ such that 
\begin{align}\label{weight-JN0}
\displaystyle\frac{1}{\displaystyle|P_r^+(z)|_a} \iint_{P_r^+(z)}|f-(f)^a_{P_r^+(z)}|^p\mu_adt
&\le C_p \big[f\big]^p_{{\rm{BMO}}(P_\frac12^+(0,1),\ \mu_adt)}, \ \forall P_r(z)\subset P_\frac12(0,1).
\end{align}
\end{lemma}

For the proof of Theorem \ref{hybrid1}, we also need the following $|y|^a$-weight version of the div-curl Lemma by Coifman-Lions-Meyer-Semmes
\cite{CLMS}, see Alvarez-Bagby-Kurtz-Perez \cite{ABKP} for the proof.

\begin{lemma}\label{hardy}
If $f\in L^2((\R^{m+1}\times \R;\mu_adt),\R^{m+1})$ is curl-free, i.e. ${\rm{curl}}_X f(X,t)=0$ for a.e. $t$,
 and $g\in L^2((\R^{m+1}\times \R;\mu_{-a}dt), \R^{m+1})$
is divergence free, i.e. ${\rm{div}}_{X}g(X,t)=0$ for a.e. $t$,  then $f\cdot g\in \mathcal{H}^1(\R^{m+1}\times \R)$, the Hardy space,  and
\begin{equation}\label{hardy-est}
\Big\|f\cdot g\Big\|_{\mathcal{H}^1(\R^{m+1}\times \R)}\le C\Big\|f\Big\|_{L^2(\R^{m+1}\times \R,\ \mu_adt)}
\Big\|g\Big\|_{L^2(\R^{m+1}\times \R,\ \mu_{-a}dt)}.
\end{equation}
\end{lemma}

\bigskip
Finally, we would like to remark on the function spaces to which the stationary heat flow of $s$-harmonic maps belongs.
Recall from Definition 1.1 that on the one hand, we assume $u\in \dot{H}^{s,\frac{s}2}_{x,t}(\R^m, \mathbb{S}^{\ell-1})$;
on the other hand, since the stationarity of $u$
is formulated by requiring that its extension map $U:\overline{\R^{m+1}_+}\times (0,\infty)\to \mathbb{R}^{\ell}$ satisfies both the $y^a$-weighted parabolic equation \eqref{eq:fhf3} and the Definition
\ref{stationary00} of stationarity condition. The validity naturally implies 
$$
U\in L^\infty\big(\R_+;  \dot{H}^1((\R^{m+1}_+,\ \mu_a), \R^{\ell})\big)\cap \dot{H}^1\big(\R_+;  L^2((\R^{m+1}_+,\ \mu_a), \R^{\ell})\big),
$$
which, with the help of the trace theorem by \cite{Nek} on $y^a$-weighted Sobolev spaces, implies that $u$ enjoys better Sobolev regularity in $t$-direction.
In particular, $u\in H^{s,s}_{(x,t), {\rm{loc}}}\big(\R^m\times (0,\infty),\mathbb{S}^{\ell-1}\big)$,  
with $(-\Delta)^{\frac{s}2} u\in L^\infty(\R_+, L^2(\R^m))$ and $(|D|_t)^s u\in L^2(\R_+, L^2(\R^m))$.

\bigskip
\noindent{\bf Acknowledgments}. 
Wang is partially supported by NSF DMS 2453789 and Simons Travel Grant TSM-00007723. Schikorra is partially supported by NSF Career DMS-2044898.

\end{document}